\documentclass[10pt]{amsart}

\usepackage{amsmath, amsfonts, amsthm, amssymb, graphicx, fullpage, enumerate}

\newtheorem{theorem}{Theorem}[section]      
\newtheorem{lemma}[theorem]{Lemma}

\newtheorem{proposition}[theorem]{Proposition}
 
\newtheorem*{main theorem}{Main Theorem}

\theoremstyle{remark}     
   
\theoremstyle{definition}  
\newtheorem{definition}[theorem]{Definition}   
      
\def\cQ{\mathcal{Q}}

\def\d{d}
\def\n{\mathbf{n}}
 
\def\i{\textnormal{i}}
\def\T{\mathbb{T}}         
\def\N{\mathbb{N}}     
          
\def\R{\mathbb{R}}     
\def\Z{\mathbb{Z}}
     
\def\C{\mathbb{C}}  
\def\P{\mathcal{P}}

\def\bar#1{\overline{#1}} 

\def\norm#1{\|#1\|}

\begin{document}
\title{A Maximal Extension of the Best-Known Bounds for the Furstenberg-S\'ark\"ozy Theorem} 
\author{Alex Rice}

\address{Department of Mathematics, Millsaps College, Jackson, MS 39210}
\email{riceaj@millsaps.edu} 
\begin{abstract}  We show that if $h\in \Z[x]$ is a polynomial of degree $k\geq 2$ such that $h(\N)$ contains a multiple of $q$ for every $q\in\N$, known as an \textit{intersective polynomial}, then any subset of $\{1,2,\dots,N\}$ with no nonzero differences of the form $h(n)$ for $n\in\N$ has density at most a constant depending on $h$ and $c$ times $(\log N)^{-c\log\log\log\log N}$, for any $c<(\log((k^2+k)/2))^{-1}$. Bounds of this type were previously known only for monomials and intersective quadratics, and this is currently the best-known bound for the original Furstenberg-S\'ark\"ozy Theorem, i.e. $h(n)=n^2$. The intersective condition is necessary to force any density decay for polynomial difference-free sets, and in that sense our result is the maximal extension of this particular quantitative estimate. Further, we show that if $g,h\in \Z[x]$ are intersective, then any set lacking nonzero differences of the form $g(m)+h(n)$ for $m,n\in \N$ has density at most $\exp(-c(\log N)^{\mu})$, where $c=c(g,h)>0$, $\mu=\mu(\deg(g),\deg(h))>0$, and $\mu(2,2)=1/2$. We also include a brief discussion of sums of three or more polynomials in the final section.
\end{abstract}
\maketitle
\setlength{\parskip}{5pt}   

\section{Introduction}  

\subsection{Background} Lov\'asz posed the following question: If $A\subseteq \N$ contains no pair of distinct elements that differ by a perfect square, must it be the case that $$\lim_{N\to \infty} \frac{|A\cap [1,N]|}{N}=0 \ ?$$ Here and throughout we use $[1,N]$ to denote $\{1,2,\dots,N\}$, and we use $|X|$ to denote the size of a finite set $X$. Furstenberg \cite{Furst} answered this question in the affirmative via ergodic theory, specifically his correspondence principle, but obtained no quantitative information on the rate at which the density must decay. Independently, S\'ark\"ozy \cite{Sark1}  showed via Fourier analysis, specifically a density increment argument driven by the Hardy-Littlewood circle method, that if $A\subseteq [1,N]$ contains no nonzero square differences, then \begin{equation}\label{sark1bound} \frac{|A|}{N} \ll \left(\frac{(\log \log N)^2}{\log N}\right)^{1/3}.\end{equation} Throughout the paper we use $\log$ to denote the natural logarithm. We use ``$\ll$" to denote ``less than a constant times", with subscripts indicating on what parameters, if any, the implied constant depends.

\subsection{Improvements and extensions} Using a more intricate Fourier analytic argument, Pintz, Steiger, and Szemer\'edi \cite{PSS} improved (\ref{sark1bound}) to \begin{equation}\label{PSSbound} \frac{|A|}{N} \ll (\log N)^{-c\log\log\log\log N}, \end{equation} with $c=1/12$. 

A natural generalization of Lov\'asz's question is to extend from perfect squares to the image of more general polynomials. Balog, Pelik\'an, Pintz, and Szemer\'edi \cite{BPPS} extended (\ref{PSSbound}) to sets with no $k$-th power differences for a fixed $k\in \N$, with $c=1/4$ and the implied constant depending on $k$. 

More generally, to hope for such a result for a given nonzero polynomial  $h \in \Z[x]$, it is clearly necessary that $h(\N)$ contains a multiple of $q$ for every $q \in \N$, as otherwise there is a set $q\N$ with positive density and no differences in the image of $h$. It follows from a theorem of Kamae and Mend\`es France \cite{KMF} that this condition is also sufficient, in a qualitative sense, and in this case we say that $h$ is an \textit{intersective polynomial}. 

Equivalently, a nonzero polynomial is intersective if it has a $p$-adic integer root for every prime $p$. Examples of intersective polynomials include any nonzero polynomial with an integer root or two rational roots with coprime denominators. However, there are also intersective polynomials with no rational roots, such as $(x^3-19)(x^2+x+1)$. 

It is a theorem of Lucier \cite{Lucier}, with minor improvements exhibited by Lyall and Magyar \cite{LM} and the author \cite{thesis}, that if $h\in \Z[x]$ is an intersective polynomial of degree $k\geq 2$ and $A\subseteq[1,N]$ has no nonzero differences in the image of $h$, then \begin{equation*}\frac{|A|}{N} \ll_h \left(\frac{\log\log N}{\log N}\right)^{1/(k-1)}. \end{equation*} Further, Hamel, Lyall,  and the author \cite{HLR} extended (\ref{PSSbound}) to all intersective polynomials of degree two, for any $c<1/\log(3)$ and the implied constant depending on $c$ and the polynomial. 

Apart from the original results of Furstenberg and S\'ark\"ozy, we have primarily alluded to the best-known results in each case, all established through versions of two Fourier analytic attacks. For  an extensive literature of intermediate and related results, as well as alternative proofs, the reader may refer to (in chronological order) \cite{Sark3}, \cite{Green}, \cite{Slip},  \cite{Lucier2}, \cite{Ruz}, \cite{LM}, \cite{lipan}, \cite{Lyall}, \cite{Rice}, and \cite{taoblog}.

\subsection{Main results} Here we adapt the Fourier analytic double iteration strategy first developed in \cite{PSS} and extend (\ref{PSSbound}) to the full collection of intersective polynomials.

\begin{theorem}\label{main1} Suppose $h\in \Z[x]$ is an intersective polynomial of degree $k\geq 2$ and $A\subseteq [1,N]$. If $$a-a'\neq h(n)$$ for all distinct pairs $a,a' \in A$ and all $n\in\N$, then \begin{equation*} \frac{|A|}{N} \ll_{h,c} (\log N)^{-c\log\log\log\log N}, \end{equation*} for any $c<\left[\log\left(\frac{k^2+k}{2}\right)\right]^{-1}$.
\end{theorem}

\noindent By more closely mimicking details of \cite{BPPS} and \cite{HLR}, one may be able to take the constant $c$ in Theorem \ref{main1} to be independent of $k$, and perhaps arbitrarily close to $1/\log 3$, a natural limit of the method. For our exposition, however, this range of values for $c$ is optimal. We discuss this further at the end of Section \ref{itn2}. 

Our crucial new ingredients are motivated and discussed in Sections \ref{sievesecNew} and \ref{sievesec}, respectively, with the necessary exponential sum estimates established in Section \ref{expest}. In addition, for the interested reader, a summary of our new exponential sum estimates is provided as a stand alone theorem in Section \ref{standalone}.   

Further, we apply the exponential sum estimates established in Section \ref{expest} in a more straightforward $L^2$ density increment, exhibiting stronger density bounds on sets free of nonzero differences that are the sum of two polynomial images. 

\begin{theorem}\label{main2} Suppose $g, h \in \Z[x]$ are nonzero intersective polynomials and $A\subseteq [1,N]$. If $$a-a'\neq g(m)+h(n)$$ for all distinct pairs $a,a'\in A$ and all $m,n \in \N$, then \begin{equation*}\frac{|A|}{N} \ll_{g,h} e^{-c(\log N)^{\mu}},  \end{equation*} where $c=c(g,h)>0$, $\mu=\mu(\deg(g),\deg(h))>0$, and $\mu(2,2)=1/2$.
\end{theorem}

\noindent \textit{Remark on generality of Theorem \ref{main2}.} We note that the necessary intersective condition makes perfect sense in a multivariable setting, and analogous results should hold for every intersective integral polynomial in several variables, not just diagonal forms. Further, there do exist intersective binary diagonal forms not covered in this theorem. For example, if $p$ is a prime congruent to $1$ modulo $90090$ that is not the sum of two integer cubes (of which there are plenty), then, since $p$ is a sum of two cubes modulo $q$ for every $q\in \N$, $x^3+y^3-p$ is an intersective polynomial in two variables that cannot be expressed as the sum of two single-variable intersective polynomials.

\subsection{Lower bounds and conjectures} For $k,N\in \N$, by fixing a prime $N^{1/k}/2 \leq p \leq N^{1/k}$ and letting $$A=\{xp: 1\leq x\leq p^{k-1}\},$$ we see that $A\subseteq [1,N]$ has no nonzero $k$-th power differences. More generally, for any polynomial $h\in \Z[x]$ of degree $k$, the greedy algorithm produces a set $A\subseteq[1,N]$ satisfying $|A|\gg_h N^{1-1/k}$ with no nonzero differences in the image of $h$.

Ruzsa \cite{Ruz2} showed that if $q\in \N$ is squarefree and $B\subseteq \Z/q\Z$ has no nonzero differences that are $k$-th powers modulo $q$, then there exists $A\subseteq [1,N]$ with no nonzero $k$-th power differences satisfying $|A|\gg N^c$, where $c=(k-1+\log |B|/\log q)/k$,  which is larger than the trivial construction with which we began this section. For $k=2$, Lewko \cite{Lewko} utilized an extensive computer search and found an example with $q=205$ and $|B|=12$, yielding $c\approx 0.7334$, which is currently the best-known lower bound for the original square-difference question.

As Ruzsa remarks, if $k=2$, then $\log|B|/\log q$ cannot exceed $1/2$ if $q$ is prime, and he conjectured this to be the case for all squarefree $q$. This indicates that $N^{3/4}$ is a limitation of Ruzsa's construction for square difference-free sets. Further, an easy Fourier analytic argument yields that if $p$ is prime and $A\subseteq (\Z/p\Z)^2=G$ has no nonzero differences of the form $(t,t^2)$, a set of forbidden differences similar in density and structure to the squares  in $[1,N]$, then $|A| \ll p^{3/2}=|G|^{3/4}$. 

These observations could potentially be viewed as evidence toward $N^{3/4}$ as roughly the true threshold for avoiding square differences, but these heuristics are rather tenuous. In particular, the same Fourier analytic argument applied to $G=(\Z/p\Z)^k$ gives an upper bound of about $|G|^{1-\frac{1}{2k}}$ for sets without differences of the form $(t,t^2,\dots,t^k)$, while Ruzsa's construction can beat this exponent for certain values of $k$. All of these questions are still massively open, and many believe these thresholds grow faster than $N^{1-\epsilon}$ for any $\epsilon>0$.

\section{Preliminaries} \label{auxpre} In this section we make some preliminary definitions and observations required to execute the Fourier analytic double iteration strategy utilized to prove Theorem \ref{main1}, as well as the more straightforward density increment utilized to prove Theorem \ref{main2}. We also provide context and motivation in Section \ref{sievesecNew} for our most notable new ingredient, a polynomial specific sieve that is defined and discussed rigorously in Section \ref{sievesec}, and we summarize our sieved exponential sum estimates in a stand alone theorem in Section \ref{standalone}.

\subsection{Fourier analysis and the circle method on $\Z/N\Z$} \label{ZNZ}
We identify  subsets of the interval $[1,N]$ with subsets of the finite group $\Z_N=\Z/N\Z$, on which we utilize a normalized discrete Fourier transform. Specifically, for a function $F: \Z_N \to \C$, we define $\widehat{F}: \Z_N \to \C$ by \begin{equation*} \widehat{F}(t) = \frac{1}{N} \sum_{x \in \Z_N} F(x)e^{-2 \pi ixt/N}. \end{equation*} 
When considering a set $A\subseteq \Z_N$, we employ a common abuse of notation by letting $A(x)$ denote the characteristic function of $A$. We analyze the Fourier analytic behavior of $A$ using the Hardy-Littlewood circle method, decomposing the nonzero frequencies into two pieces: the points $t\in \Z_N$ such that $t/N$ is close to a rational with small denominator, and the complement.

\begin{definition}\label{arcs} Given $N\in \N$ and $K,Q>0$, we define, for each $q\in \N$ and $a \in [1,q]$, 
\begin{equation*} \mathbf{M}_{a,q}(N,K)=  \left\{t \in \Z_N: \left|\frac{t}{N}-\frac{a}{q}\right| < \frac{K}{N} \right\} \text{ \ and \ } \mathbf{M}_q(N,K) = \bigcup_{(a,q)=1}\mathbf{M}_{a,q}(N,K)\setminus \{0\}. \end{equation*} 
We then define $\mathfrak{M}(N,K,Q)$, the \emph{major arcs}, by \begin{equation*} \mathfrak{M}(N,K,Q) = \bigcup_{q=1}^{Q} \mathbf{M}_q(N,K), \end{equation*} and $\mathfrak{m}(N,K,Q)$, the \emph{minor arcs}, by $\mathfrak{m}(N,K,Q) = \Z_N \setminus (\mathfrak{M}(N,K,Q) \cup \{0\})$. It is important to note that as long as $2KQ^2<N$, we have that $\mathbf{M}_{a,q}(N,K) \cap \mathbf{M}_{b,r}(N,K) = \emptyset$ whenever $a/q \neq b/r$, $q,r \leq Q$.
\end{definition} 


\subsection{Auxiliary polynomials}\label{auxdef1} Suppose $h\in \Z[x]$ is an intersective polynomial. For each prime $p$, we fix a $p$-adic integer $z_p$ with $h(z_p)=0$. The objects defined below certainly depend on $h$, as well as on the choice of $p$-adic integer roots, though any choice works equally well for our purposes, and we suppress all of this dependence in the coming notation. 

 By reducing modulo prime powers and applying the Chinese Remainder Theorem, the choices of $z_p$  determine, for each natural number $d$, a unique integer $r_d \in (-d,0]$, which consequently satisfies $d \mid h(r_d)$. We define the function $\lambda$ on $\N$ by  letting $\lambda(p)=p^{m}$ for each prime $p$, where $m$ is the multiplicity of $z_p$ as a root of $h$, and then extending it to be completely multiplicative. For each $d\in \N$, we define the \textit{auxiliary polynomial}, $h_d$, by 
\begin{equation*} h_d(x)=h(r_d + dx)/\lambda(d). \end{equation*} If $p^j \mid d$ for $p$ prime and $j \in \N$, then since $r_d \equiv z_p$ mod $p^j$, we see by factoring $h$ over the $p$-adic integers that all the coefficients of $h(r_d+dx)$ are divisible by $p^{jm}$, hence each auxiliary polynomial has integer coefficients.

It is important to note that the leading coefficients of the auxiliary polynomials grow at least as quickly, up to a constant depending only on $h$, as the other coefficients. Specifically, if $h(x)=a_kx^k+\cdots+a_1x+a_0 \in \Z[x]$ with $a_k>0$, then the leading coefficient of $h_d(x)$ is $d^ka_k/\lambda(d)$, while the magnitude of every coefficient is (quite brutally) bounded by $2^kd^k(\max_{0\leq i \leq k} |a_i|)/\lambda(d).$ This fact combines usefully with the following proposition, which carefully keeps track of the extent to which a polynomial behaves like its leading term with regard to preimages.

\begin{proposition} Suppose $h(x)=a_kx^k+\cdots+a_1x+a_0 \in \Z[x]$, $a_k>0$, and let $R=(|a_0|+\cdots+|a_{k-1}|)/a_k$. Then, for any $x>0$, $$ \left|\left\{n\in \N: 0<h(n)<x\right\} \ \triangle \ [1,(x/a_k)^{1/k}]\right| \leq 3 \lfloor R \rfloor +2, $$ where $\Delta$ denotes the symmetric difference. \

\end{proposition}

\begin{proof} Fixing $h\in \Z[x]$ and $x>0$ as in the proposition, and letting $y=(x/a_k)^{1/k}$, we see that elements of the symmetric difference at hand arise in three ways: $n\in \N$ such that $h(n)\leq 0$, $y+n \in \N$ with $n>0$ such that $h(y+n)<x$, and $y-n\in \N$ with $0\leq n<y$ such that $h(y-n)\geq x$. Defining $R$ as in the proposition, it suffices to show that these sets have size at most $\lfloor R \rfloor$, $\lfloor R \rfloor+1$, and $\lfloor R \rfloor+1$, respectively. To this end, we observe that for $n\in \N$ we have $a_kn^k-Ra_kn^{k-1} \leq h(n) \leq a_kn^k+Ra_kn^{k-1}.$ In particular, $h(n)>0$ if $n>R$. Carefully expanding  $a_k(y+n)^k-Ra_k(y+n)^{k-1}$, yields $$a_ky^k+a_k \sum_{i=0}^{k-1} n^i({k \choose i+1}n-R{k-1 \choose i})y^{k-1-i} \geq x+a_k\sum_{i=0}^{k-1}{k \choose i+1} n^i(n-R)y^{k-1-i}, $$ which is at least $x$ provided that $n\geq R$. Analogously, expanding $a_k(y-n)^k+Ra_k(y-n)^{k-1}$ yields $$x +a_k  \left(R \sum_{i=0}^{k-1}{k-1 \choose i} (-n)^iy^{k-1-i}-n\sum_{i=0}^{k-1} {k \choose i+1} (-n)^iy^{k-1-i}\right). $$ We see that for $0<n<y$, both summations are alternating such that each negative term can be paired off with a positive term that is at least as large. Further, the latter summation is at least as large in magnitude as the former, hence $h(y-n)<x$ for $R<n<y$, and the proposition follows.
\end{proof}

In particular, if $\deg(h)=k$ and $b_d>0$ is the leading coefficient of $h_d$, then for any $x>0$ we have
\begin{equation} \label{symBIG}
\left|\left\{n\in \N: 0<h_d(n)<x\right\} \ \triangle \ [1,(x/b_d)^{1/k}]\right| \ll_{h} 1,
\end{equation} which will serve as a convenient observation for our exposition.   

\subsection{Inheritance proposition} We define these auxiliary polynomials to keep track of an inherited lack of prescribed differences at each step of a density increment iteration. To this end, for nonzero $h\in\Z[x]$, we define $I(h)$ to be the positive elements of $h(\N)$ if $h\in \Z[x]$ has positive leading coefficient and the negative elements of $h(\N)$ if $h\in \Z[x]$ has negative leading coefficient. The following proposition makes the aforementioned inheritance precise.

\begin{proposition} \label{inh} Suppose $h\in\Z[x]$ is intersective with positive leading coefficient, $d,q\in \N$, and $A\subseteq \N$.

\noindent If $(A-A)\cap I(h_d)=\emptyset$ and $A'\subseteq \{a: x+\lambda(q)a \in A\}$, then $(A'-A')\cap I(h_{qd})=\emptyset$.
\end{proposition}
\begin{proof} Suppose that $A\subseteq \N$, $A'\subseteq \{a : x+\lambda(q)a \in A\}$, and $$a-a'=h_{qd}(n)=h(r_{qd}+qdn)/\lambda(qd)>0$$ for some $n\in \N$, $a, a' \in A'$. By construction we see that $r_{qd}\equiv r_d$ mod $d$, so there exists $s\in \Z$ such that $r_{qd}=r_d+ds$, and therefore $$0<h_d(s+qn)=\frac{h(r_d+d(s+qn))}{\lambda(d)}=\lambda(q)h_{qd}(n)=\lambda(q)a-\lambda(q)a' \in A-A,$$ hence $(A-A)\cap I(h_d) \neq \emptyset$, and the contrapositive is established. 
\end{proof} 

\noindent This single polynomial inheritance proposition is all we need for the purposes of proving Theorem \ref{main1}, but we require a nominally generalized version of Proposition \ref{inh} for results involving multiple polynomials, which we include in Section \ref{singsec}.

\subsection{Sieve motivation}\label{sievesecNew} In this context, the guiding principle of the Hardy-Littlewood circle method is that if $h\in \Z[x]$, then the \textit{Weyl sum} \begin{equation} \label{weyls} \sum_{n=1}^Me^{2\pi \textnormal{i}h(n)\alpha}\end{equation} is much smaller than the trivial bound $M$, unless $\alpha$ is well-approximated by a rational number with small denominator. 

On a coarse scale, this principle is captured by combining the pigeonhole principle with Weyl's inequality (see Lemma \ref{weyl2}), but for a more refined treatment, we must address the following question: If $\alpha$ \textit{is} close to a rational with a small denominator, for example a denominator smaller than a tiny power of $M$, can we beat the trivial bound at all? 

This question turns out to be quite straightforward, as under these conditions (\ref{weyls}) has a convenient asymptotic formula, and the gain from the trivial bound resides in a local version of the sum. Specifically, if $\alpha$ is close to $a/q$ with  $q$ small, then, up to a small error, the magnitude of (\ref{weyls}) is at most $M$ times  \begin{equation}\label{gsintro} q^{-1}\sum_{s=0}^{q-1}e^{2\pi \text{i}h(s)a/q}. \end{equation} 

In general, the best decay we can hope for from (\ref{gsintro}) is $q^{-1/k}$, where $k=\deg(h)$ (see \cite{Chen} for example). Unfortunately, it is completely vital to the double iteration method developed in \cite{PSS} to establish (\ref{PSSbound}) that one has decay at or near $q^{-1/2}$ (what we refer to as \textit{square root cancellation}) for small denominators $q$. Of course, this is not an obstacle if $k=2$, a fact exploited in \cite{HLR} to establish (\ref{PSSbound}) for all intersective quadratics.

In \cite{BPPS}, Balog, Pelik\'an, Pintz, and Szemer\'edi observed that by sieving the initial set of inputs, letting $W$ equal a product of small primes and considering only inputs coprime to $W$, the gain from the trivial bound for $\alpha$ near $a/q$ with small $q$ is given roughly by \begin{equation}\label{gsintro2} \phi(q)^{-1}\sum_{\substack{s=0 \\ (s,q)=1}}^{q-1}e^{2\pi \text{i}h(s)a/q}. \end{equation} In the case of $h(n)=n^k$, one then exploits that, for all but the primes dividing $k$, units modulo $p^j$ for $j\geq 2$ are $k$-th power residues if and only if they reduce to $k$-th power residues modulo $p$. This allows one to use orthogonality of characters and reduce to the case of squarefree $q$, where there is always square root cancellation.

However, a generic intersective polynomial need not have the root-lifting properties necessary to establish square root cancellation in (\ref{gsintro2}), and so sieving away small prime factors is not a winning strategy for maximally extending (\ref{PSSbound}). Alternatively, given an intersective polynomial $h\in \Z[x]$, we utilize a non-traditional sieve dependent on $h$, essentially defined by removing the zeros of the derivative $h'$ modulo small primes, in order to mimic the aforementioned lifting observation using Hensel's Lemma. 

\subsection{Sieve definitions and observations} \label{sievesec} For an intersective polynomial $h\in\Z[x]$ and each prime $p$ and $d\in \N$, we define $\gamma_{d}(p)$ to be the smallest power such that the derivative $h_d'$ is not identically zero modulo $p^{\gamma_{d}(p)}$, and we let $j_d(p)$ denote the number of roots of $h_d'$ modulo $p^{\gamma_{d}(p)}$.

\noindent Then, for $d\in \N$ and $Y>0$ we define $$W_d(Y)=\left\{n\in \N: h_d'(n) \not\equiv 0 \text{ mod } p^{\gamma_{d}(p)} \text{ for all } p\leq Y \right\}.$$
  
\noindent \noindent In the absence of a subscript $d$ in the usage of $\gamma(p), j(p),$ and $W(Y)$, we assume $d=1$, in which case the definitions make sense even for non-intersective polynomials. Further, for any $g\in \Z[x]$ and $q\in \N$, we define $$W^{q}(Y)=\left\{n\in \N: g'(n) \not\equiv 0 \text{ mod }p^{\gamma(p)} \text{ for all } p\leq Y, \ p^{\gamma(p)}\mid q \right\}.$$  The size of $W(Y)$ can be estimated with a standard Brun sieve calculation, which for completeness we include below.\\

\begin{proposition}\label{brunprop} If $g\in \Z[x]$ and $X,Y>0$ with $c\log X \geq \log Y \log \log Y$, then \begin{equation}\label{sieve} \left|[1,X]\cap W(Y)\right| = X\prod_{p\leq Y} \left(1-\frac{j(p)}{p^{\gamma(p)}} \right)+O\left(Xe^{-c\frac{\log X}{\log Y}} \right), \end{equation} where $c>0$ depends only on $\deg(g)$ and the collection of moduli for which $g'$ is identically zero.

\end{proposition}

\begin{proof} We fix $g\in \Z[x]$ and $X>0$, and for each collection $p_1<p_2<\cdots<p_s$ of primes, we let \begin{equation*}\mathcal{A}_{p_1\cdots p_s}(X) = \left|\left\{1\leq n \leq X: g'(n) \equiv 0\text{ mod } p_i^{\gamma(p_i)} \ \text{for all } 1\leq i \leq s  \right\}\right|.\end{equation*} 

\noindent Fixing $Y>0$ and letting $r$ denote the number of primes that are at most $Y$, we have by the Chinese Remainder Theorem and the inclusion-exclusion principle that \begin{equation}\label{brunalt1} \left|[1,X] \cap  W(Y)\right|= \sum_{s=0}^r  (-1)^s \sum_{p_1<\dots < p_s} \mathcal{A}_{p_1\cdots p_s}(X),\end{equation} and moreover the true size lies between any two consecutive truncated alternating sums in $s$. Further, we know that \begin{equation}\label{AP1} \mathcal{A}_p = \frac{j(p)X}{p^{\gamma(p)}} + R_p, \end{equation} where $|R_p|\leq j(p)$. \\

\noindent We now observe that $j(p)\leq k-1$ if $\gamma(p)=1$, and trivially $j(p) \leq p^{\gamma(p)}$. In particular, $j(p)$ can be bounded above in terms of only $\deg(g)$ and the collection of moduli for which $g'$ is identically zero, which allows us to apply the Chinese Remainder Theorem again and extend (\ref{AP1}) to \begin{equation} \label{AP2}  \mathcal{A}_{p_1\cdots p_s} = X \prod_{i=1}^s\frac{j(p_i)}{p_i^{\gamma(p_i)}} + R_{p_1\cdots p_s}, \end{equation} where $|R_{p_1\cdots p_s}| \leq C^s$.\\

\noindent For this and the remainder of the proof we let $C$ denote a positive constant, which may change from line to line but depends only on $\deg(g)$ and the collection of moduli for which $g'$ is identically zero.  

\noindent For any even $0\leq t \leq r$, we have by (\ref{brunalt1}) and (\ref{AP2}) that 
 
\begin{align*}  \left|[1,X] \cap  W(Y)\right| &\leq \sum_{s=0}^t  (-1)^s \sum_{p_1<\dots < p_s} \mathcal{A}_{p_1\cdots p_s}(X) \\\\ & = \sum_{s=0}^t  (-1)^s \sum_{p_1<\dots < p_s \leq Y}  X \prod_{i=1}^s\frac{j(p_i)}{p_i^{\gamma(p_i)}} + R_{p_1\cdots p_s} \\\\ & = X\prod_{p\leq Y} \left(1-\frac{j(p)}{p^{\gamma(p)}} \right) - X\sum_{s=t+1}^r  (-1)^s \sum_{p_1<\dots < p_s \leq Y}   \prod_{i=1}^s\frac{j(p_i)}{p_i^{\gamma(p_i)}} + \sum_{s=0}^t \sum_{p_1<\dots < p_s \leq Y} R_{p_1\cdots p_s} \\\\ &\leq X\prod_{p\leq Y} \left(1-\frac{j(p)}{p^{\gamma(p)}} \right) + X \sum_{s=t+1}^r  \sum_{p_1<\dots < p_s \leq Y}   \prod_{i=1}^s\frac{C}{p_i^{\gamma(p_i)}} + \sum_{s=0}^t { r \choose s } C^s \\\\ &= X\prod_{p\leq Y} \left(1-\frac{j(p)}{p^{\gamma(p)}} \right) + E_1 + E_2,
\end{align*}
and we obtain analogous lower bounds by choosing odd $0\leq t \leq r$. To control $E_1$, we observe that \begin{equation*} \sum_{p_1<\dots < p_s \leq Y}   \prod_{i=1}^s\frac{C}{p_i^{\gamma(p_i)}} \leq \frac{1}{s!} \left(\sum_{p\leq Y} \frac{C}{p_i}\right)^s.\end{equation*} 
Using the standard fact that $\sum_{p\leq Y} 1/p \ll \log\log Y$, we then have \begin{equation*} E_1 \leq X \sum_{s>t} \frac{(C\log \log Y)^s}{s!}. \end{equation*} If $t>2C\log\log Y$, then each term of this series is at most half the previous, so the full tail is at most twice the first term. We then use the bound $t! \geq (t/e)^t$ to establish \begin{equation} \label{e1est} E_1 \leq X \left(\frac{C\log \log Y}{t}\right)^t. \end{equation}

\noindent To control $E_2$, we use the trivial bound ${ r \choose s } \leq r^s/s!$ to see that \begin{equation*} E_2\leq \sum_{s=0}^t \frac{(Cr)^s}{s!}. \end{equation*} Since $t\leq r$, each term of this series is at least twice the previous, so the full sum is bounded by double the final term. Since $r\ll Y/\log Y$ by the Prime Number Theorem, and again $t!\geq (t/e)^t$, we have \begin{equation} \label{e2est} E_2 \leq \left(\frac{CY}{t\log Y}\right)^t. \end{equation} We now finish the proof by making a good choice for $t$. If $Y\leq \log X$, then we can choose $t=r$, in which case $E_1=0$ and $E_2 \leq C^r \leq \exp(O(\log X/ \log \log X))$, which satisfies the proposition with room to spare. Finally, if $Y> \log X$, then we choose $t=\frac{\log X}{C\log Y}.$ By our hypotheses on $X$ and $Y$, this choice of $t$ satisfies the lower bound needed for (\ref{e1est}), and substituting this choice into (\ref{e1est}) and (\ref{e2est}) yields the desired error bounds. 
\end{proof} 

As it applies to auxiliary polynomials, the conclusion of Proposition \ref{brunprop}, as well as other steps in our future arguments, depends on the collection of moduli for which $h_d'$ is identically zero. The following observations assure that this collection of moduli remains under control. 
 
\begin{proposition} \label{idzero} If $g(x)=a_0+a_1x+\cdots+a_kx^k \in \Z[x]$ is identically zero modulo $q\in \N$, then $$q \mid k!\gcd(a_0,\dots,a_k). $$
\end{proposition}  

\begin{proof} We first note that $g$ is identically zero modulo $q$ if and only if the polynomial $g/q$ is integer-valued. In this case, since the binomial coefficients $$\binom{x}{j}=\frac{x(x-1)\dots(x-j+1)}{j!} $$ form a $\Z$-basis for integer-valued polynomials, we can write $$g(x)=\sum_{j=0}^k qb_j\binom{x}{j}$$ for $b_0,\dots,b_k \in \Z$. In particular, by clearing denominators we see that the coefficients of $k!g$ are all divisible by $q$, and the proposition follows.
\end{proof}

For $h(x)=a_0+a_1x+\cdots+a_kx^k\in \Z[x]$, we define the \textit{content} of $h$ by $$\text{cont}(h)=\gcd(a_1,\dots,a_k), $$ noting that any common factor of the coefficients of $h'$ divides $k!\text{cont}(h)$.

The last hurdle in controlling the set of ``bad moduli" in our sieve, as well as a major issue in our future exponential sum estimates, is the possibility that the coefficients of the auxiliary polynomials $h_d$ gain larger and larger common factors as $d$ grows, but the following lemma due to Lucier asserts that this is not case.\\

\begin{lemma}[Lemma 28, \cite{Lucier}] \label{content} If $h\in \Z[x]$ is intersective with $\deg(h)=k$, then for every $d\in \N$, \begin{equation*} \textnormal{cont}(h_d) \leq |\Delta(h)|^{(k-1)/2}\textnormal{cont}(h), \end{equation*} where $$\Delta(h)=a^{2k-2}\prod_{i\neq i'} (\alpha_i-\alpha_{i'})^{e_ie_{i'}}$$ if $h$ factors over the complex numbers as $$h(x)=a(x-\alpha_1)^{e_1}\dots(x-\alpha_r)^{e_r}$$ with all the $\alpha_i$'s distinct.
\end{lemma}

For intersective $h\in\Z[x]$ with $\deg(h)=k$, we have now established control over not only the error term in the size of $W_d(Y)$, but also the main term, since Proposition \ref{idzero}, Lemma \ref{content}, and the fact that $h_d'$ has at most $k-1$ roots modulo every prime at which it is not identically zero,  give \begin{equation}\label{logk} \prod_{p\leq Y} \left(1-\frac{j_d(p)}{p^{\gamma_d(p)}} \right) \gg_h \prod_{k\leq p\leq Y} \left(1-\frac{k-1}{p} \right) \gg (\log Y)^{1-k}\end{equation} for all $d \in \N$ and $Y\geq 2$.

\subsection{Summary of new exponential sum estimates} \label{standalone} In Section \ref{expest}, we combine new and old techniques to establish the sieved exponential sum estimates necessary to prove Theorem \ref{main1} and additional results. These estimates are obtained through a sequence of lemmas presented in the context of the larger proof, so we separately present a summary here in case the estimates are of independent interest to the reader.
 
\noindent For the following theorem, we utilize all the sieve-related notation and definitions from Section \ref{sievesec}. 

\begin{theorem} For $k \geq 2$, $g(x)=a_0+a_1x+\cdots+a_kx^k \in \Z[x]$, and $a,q\in \N$, the following estimates hold:

\begin{enumerate}[(i)]  \item \textnormal{\textbf{Major arc estimate:}} If $X,Y\geq 2$ with $c\log (X/q) \geq \log Y \log \log Y$, $J=|a_0|+|a_1|+\cdots+|a_k|$, and $\alpha=a/q+\beta$, then \begin{align*}\sum_{\substack{n=1 \\ n \in W(Y)}}^Xg'(n)e^{2\pi \textnormal{i}g(n)\alpha}&=\frac{1}{q} \prod_{\substack{ p\leq Y \\ p^{\gamma(p)}\nmid q}}\left(1-\frac{j(p)}{p^{\gamma(p)}}
\right)\sum_{\substack{s=0 \\ s \in W^{q}(Y)}}^{q-1}e^{2\pi \textnormal{i} g(s)\frac{a}{q}}\int_0^Xg'(x)e^{2\pi\textnormal{i}g(x)\beta}dx\\\\& + O\left(kJX^ke^{-c\frac{\log\left(\frac{X}{q}\right)}{\log Y}}(1+JX^k|\beta|)\right),\end{align*} where $c=c(k,\textnormal{cont(g)})>0$. \\ \item \textnormal{\textbf{Square root cancellation:}} If $(a,q)=1$ and $Y>0$, then \begin{equation*}\left|\sum_{\substack{s=0 \\ s \in W^{q}(Y)}}^{q-1}e^{2\pi \textnormal{i} g(s)a/q}\right| \ll_k \gcd(\textnormal{cont}(g),q)^{3} C^{\omega(q)}\begin{cases} q^{1/2} &\text{if }q\leq Y \\ q^{1-1/k} &\text{else} \end{cases}, \end{equation*} where $C=C(k)$ and $\omega(q)$ is the number of distinct prime factors of $q$. \\ \item \textnormal{\textbf{Minor arc estimate:}} If $a_k>0$, $X,Y,Z \geq 2$, $(a,q)=1$ and $|\alpha-a/q|<q^{-2}$, then $$\left|\sum_{\substack{n=1 \\ n \in W(Y)}}^X e^{2\pi \i g(n)\alpha} \right| \ll_{k} \textnormal{cont}(g)^5(\log Y)^{ek} X\left(e^{-\frac{\log Z}{\log Y}}+\left(a_{k}\log^{k^2}(a_{k}qX)\left(q^{-1}+\frac{Z}{X}+\frac{qZ^k}{a_kX^k}\right) \right)^{2^{-k}} \right). $$
\end{enumerate}

\end{theorem}

\noindent These estimates are obtained mostly by combining elementary sieve methods with traditional circle method techniques. The main novelty, the exponent $1/2$ in estimate $(ii)$, is established using Hensel's Lemma.

\section{The Double Iteration Method: Proof of Theorem \ref{main1}} \label{mainsec} In this section we execute an adapted version of the double iteration method originally developed in \cite{PSS} and prove Theorem \ref{main1}. We begin with an outline.

\subsection{Overview of the argument}   \label{over}
We initially observe that if $h\in \Z[x]$ is an intersective polynomial, which by symmetry of difference sets we can assume has positive leading coefficient, and $$(A-A)\cap I(h)=\emptyset$$ for a set $A\subseteq [1,N]$, then we can apply the circle method to show that this unexpected behavior implies substantial $L^2$  mass of $\widehat{A}$ over nonzero frequencies near rationals with small denominator. 

At this point, the traditional method, which we employ in Section \ref{singsec} to prove Theorem \ref{main2}, is to use the pigeonhole principle to conclude that there is one single denominator $q$ such that $\widehat{A}$ has $L^2$ concentration around  rationals with denominator $q$. From this information, one can conclude that $A$ has increased density on a long arithmetic progression with step size an appropriate multiple of $q$, leading to a new denser set with an inherited lack of polynomial differences  and continued iteration. 

Pintz, Steiger, and Szemer\'edi \cite{PSS}  observed that pigeonholing to obtain a single denominator $q$ is a potentially wasteful step. We follow their approach, observing the following dichotomy:  
\begin{quotation}
\noindent \textbf{Case 1.} There is a single denominator $q$ such that $\widehat{A}$ has extremely high $L^2$ concentration, greater than yielded by the pigeonhole principle, around rationals with denominator $q$. This leads to a very large density increment on a long arithmetic progression. 
\medskip
 
\noindent \textbf{Case 2.} The $L^2$ mass of $\widehat{A}$ on the major arcs is spread over many denominators. In this case, an iteration procedure using the ``combinatorics of rational numbers" can be employed to build a large collection of frequencies at which $\widehat{A}$ is large, then Plancherel's identity is applied to bound the density of $A$.
\end{quotation}

\noindent Philosophically, Case 1 provides more structural information about the original set $A$ than Case 2 does. The downside is that the density increment procedure yields a new set and potentially  a new polynomial, while the iteration in Case 2 leaves these objects fixed. 
With these cases in mind, we can now outline the argument, separated into two distinct phases. 
\begin{quotation}
\noindent \textbf{Phase 1} (The Outer Iteration): Given a set $A$ and an intersective polynomial $h\in \Z[x]$ with $(A-A)\cap I(h)= \emptyset,$ we ask if the set falls into Case 1 or Case 2 described above. If it falls into Case 2, then we proceed to Phase 2. 
\medskip

If the set $A$ falls into Case 1, then the density increment procedure yields a small $q\in \N$ and a new subset $A_1$ of a slightly smaller interval with significantly greater density, and  $$(A_1-A_1) \cap I(h_q)=\emptyset.$$
We can then iterate this process as long as the resulting interval is not too small, and the dichotomy holds as long as the coefficients of the auxiliary polynomial are not too large. We show that if the resulting sets remain in Case 1 until the interval shrinks down or the coefficients grow to the limit, then the density of the original set $A$ must have satisfied a bound stronger than the one purported in Theorem \ref{main1}. Contrapositively, we assume that the original density does not satisfy this stricter bound, and we conclude that one of the sets yielded by the density increment procedure must lie in a large interval, have no differences in $I(h_d)$ for reasonably small $d$, and fall into Case 2. We call that set $B\subseteq[1,L]$. 
\medskip
\end{quotation}

\noindent We now have a set $B \subseteq [1,L]$ with $(B-B)\cap I(h_d)=\emptyset$ which falls into Case 2, so we can adapt the strategy of \cite{PSS}, \cite{BPPS}, and \cite{HLR}. It is in this phase that we use the sieve outlined in Section \ref{sievesec}, as the method breaks down without square root cancellation on the major arcs.

\begin{quotation} \textbf{Phase 2} (The Inner Iteration): We prove that given a frequency $s\in \Z_L$ with $s/L$ close to a rational $a/q$ such that $\widehat{B}(s)$ is large, there are lots of nonzero frequencies $t\in \Z_L$ with $t/L$ close to rationals $b/r$ such that $\widehat{B}(s+t)$ is almost as large. This intuitively indicates that a set $P$ of frequencies associated with large Fourier coefficients can be blown up to a much larger set $P'$ of frequencies associated with nearly as large Fourier coefficients. The only obstruction to this intuition is the possibility that there are many pairs $(a/q,b/r)$ and $(a'/q',b'/r')$ with $$\frac{a}{q}+\frac{b}{r} = \frac{a'}{q'}+\frac{b'}{r'}.$$ Observations made in \cite{PSS} and \cite{BPPS} on the combinatorics of rational numbers demonstrate that this potentially harmful phenomenon can not occur terribly often. 
\medskip

Starting with the trivially large Fourier coefficient at $0$, this process is applied as long as certain parameters are not too large, and the number of iterations is ultimately limited by the growth of the divisor function. Once the iteration is exhausted, we use the resulting set of large Fourier coefficients and Plancherel's Identity to get the upper bound on the density of $B$, which is by construction larger than the density of the original set $A$, claimed in Theorem \ref{main1}.
\end{quotation}

\subsection{Reduction to two key lemmas}  \label{twolemmas}

For the remainder of Section \ref{mainsec}, we fix an intersective polynomial $h \in \Z[x]$ with positive leading coefficient and $\deg(h)=k\geq 2$, and we let $$\rho=2^{-10k}.$$ We fix $\epsilon>0$ and let $$m=k^2+k+2\epsilon,$$ with the goal of establishing Theorem \ref{main1} with $c=(1-17\epsilon)/\log(m/2)$. We note that any dependence of implied constants on $m$ or is actually just dependence on $k$ and $\epsilon$. 

We also fix a natural number $N$, and we let $$\cQ=(\log N)^{\epsilon \log\log\log N}.$$ We preemptively assume that $N$ is sufficiently large at all junctures with respect to $h$ and $\epsilon$, an assumption that is permitted because the implied constant in Theorem \ref{main1} is allowed to depend on those parameters.  We deduce Theorem \ref{main1} from two key lemmas, corresponding to the two phases outlined in the overview in Section \ref{over}, the first of which yields a set with substantial Fourier $L^2$ mass distributed over rationals with many small denominators.

\begin{lemma} \label{outer} Suppose $A\subseteq [1,N]$ with $|A|=\delta N$ and $(A-A)\cap I(h)= \emptyset$. If \begin{equation} \label{behrend} \delta \geq e^{-(\log N)^{\epsilon/4}},\end{equation} then there exists  $B \subseteq [1,L]$ and $d \leq N^{\rho/2}$ satisfying $ N^{1-\rho} \leq L \leq N$, $|B|/L=\sigma \geq \delta$, and $$(B-B) \cap I(h_d)= \emptyset.$$ Further, $B$ satisfies $|B\cap[1,L/2]|\geq \sigma L/3$ and  \begin{equation} \label{maxmass} \max_{q\leq \cQ}\sum_{t\in \mathbf{M}_q(L,\cQ)}  |\widehat{B}(t)|^2 \leq \sigma^2(\log N)^{-1+\epsilon}. \end{equation} 
\end{lemma}  

The second lemma corresponds to the iteration scheme in which a set of large Fourier coefficients from distinct major arcs is blown up in such a way that the relative growth of the size of the set is much greater than the relative loss of pointwise mass.

\begin{lemma} \label{itn} Suppose $B\subseteq [1,L]$ and $d\in \N$ are as in the conclusion of Lemma \ref{outer}, let $B_1=B\cap[1,L/2]$,  and suppose $\sigma \geq \cQ^{-1/m}$. Given $U,V,K \in \N$ with $\max \{U,V,K\}\leq \cQ^{1/m}$ and a set \begin{equation*} P \subseteq \left\{t \in \bigcup_{q=1}^{V} \mathbf{M}_q(L,K)\cup\{0\} : |\widehat{B_1}(t)|  \geq \frac{\sigma}{U} \right\} \end{equation*} satisfying \begin{equation} \label{dis}| P \cap \mathbf{M}_{a,q}(L,K)| \leq 1 \text{ \ whenever \ } q \leq V,\end{equation} there exist $U',V',K' \in \N$ with $\max\{U',V',K'\} \ll_{h,\epsilon} (\max \{U,V,K\})^{m/2}\sigma^{-(m-1)/2}$ and a set  \begin{equation} \label{P'1} P' \subseteq \left\{t \in \bigcup_{q=1}^{V'} \mathbf{M}_q(L,K')\cup\{0\} : |\widehat{B_1}(t)|  \geq \frac{\sigma}{U'} \right\}\end{equation} satisfying \begin{equation} \label{P'2} | P' \cap \mathbf{M}_{a,q}(L,K')| \leq 1 \text{ \ whenever \ } q\leq V'\end{equation} and \begin{equation} \label{P'3} \frac{|P'|}{(U')^2} \geq \frac{|P|}{U^2}(\log N)^{1-16\epsilon}.\end{equation}
\end{lemma}

\subsection{Proof of Theorem \ref{main1}} \label{hardproof} In order to establish Theorem \ref{main1}, we can assume that \begin{equation*}\delta \geq (\log N)^{-\log\log\log\log N}. \end{equation*} Therefore, Lemma \ref{outer} produces a set $B$ of density $\sigma \geq \delta$ with the stipulated properties, and we set $P_0=\{0\}$, $U_0=3$, and $V_0=K_0=1$. Then, Lemma \ref{itn} yields, for each $n$, a set $P_n$ with parameters $U_n,V_n,K_n$ such that \begin{equation*}\max\{U_n,V_n,K_n\} \leq (C\log N)^{(m/2)^{n+1}\log\log\log\log N}, \end{equation*} where $C=C(h,\epsilon)$, and \begin{equation*}\frac{1}{\sigma} \geq \frac{1}{\sigma^2}\sum_{t\in P_n} |\widehat{B_1}(t)|^2 \geq \frac{|P_n|}{U_n^2} \gg (\log N)^{n(1-16\epsilon)}, \end{equation*} where the left-hand inequality comes from Plancherel's Identity, as long as $\max\{U_n,V_n,K_n\} \leq \cQ^{1/m}$. This holds with $n=(1-\epsilon)(\log\log\log\log N) /\log (m/2)$, as $(m/2)^{n+1} \leq (\log\log\log N)^{1-\epsilon/2}$, and the theorem follows. 
\qed

\subsection{The Outer Iteration} We begin the first phase with the following standard $L^2$ density increment lemma, and we provide a proof for the sake of completeness.

\begin{lemma} \label{dinczn} Suppose $B\subseteq [1,L]$ with $|B|=\sigma L$. If $0<\theta \leq 1$, $q\in\N$, $K>0$, and
$$\sum_{t\in \mathbf{M}'_q(L,K)}|\widehat{B}(t)|^2 \geq \theta\sigma^2,$$ then there exists an arithmetic progression \begin{equation*}P=\{x+\ell q : 1\leq \ell \leq L'\}
\end{equation*} with $qL' \gg \sigma\min\{\theta, K^{-1}\}L$ and $|B\cap P|\geq \sigma(1+\theta/16)L'$. 
\end{lemma}

\begin{proof}Suppose $B\subseteq [1,L]$ with $|B|=\sigma L$. Suppose further that  
\begin{equation}\label{qmass} \sum_{t\in \mathbf{M}'_q(L,K)}|\widehat{B}(t)|^2  \geq \theta\sigma^2 
\end{equation} 
and let $P = \{q,2q, \dots, Xq\} \subseteq \Z_L$ with $X= \lfloor\min\{\theta,  K^{-1}\} L/16 q\rfloor$. We will show that either all or a portion of some translate of $P$, lifted to the integers, satisfies the conclusion of Lemma \ref{dinczn}, with either $L'=X$ or $L'=\lceil \sigma X/2 \rceil$. We note that for $t \in \Z_L$,
\begin{equation}\label{Phat} L|\widehat{P}(t)|=\Big| \sum_{\ell=1}^{X} e^{-2\pi i \ell q t/L} \Big| \geq  X- \sum_{\ell=1}^X |1-e^{-2\pi i \ell q t/L}| \geq  X-2\pi X^2 \norm{qt/L},
\end{equation} where $\norm{\cdot}$ denotes the distance to the nearest integer. Further, if $t \in \mathbf{M}'_q(L,K)$, then \begin{equation}\label{dni} \norm{qt/L} \leq qK/L \leq 1/4\pi X.
\end{equation} Therefore, by (\ref{Phat}) and (\ref{dni}) we have 
\begin{equation}\label{L/2} L|\widehat{P}(t)| \geq X/2 \quad \text{for all} \quad t \in \mathbf{M}_q(L,K).
\end{equation}
To investigate the bias of $B$ toward translations of $P$, we introduce the balanced function $f_B:\Z_L \to \R$ defined by $f_B(x)=B(x)-\sigma$, noting by orthogonality of characters that \begin{equation} \label{bfo} \widehat{f_B}(t) = \begin{cases} 0 & \text{if } t=0 \\ \widehat{B}(t) & \text{else} \end{cases}. \end{equation}By (\ref{qmass}) and (\ref{L/2})  we see
\begin{equation}\label{convmass1}  \theta\sigma^2 \leq \sum_{t\in \mathbf{M}'_q(L,K)}|\widehat{B}(t)|^2 \leq \frac{4L^2}{X^2} \sum_{t\in \Z_L \setminus \{0\}}|\widehat{B}(t)|^2|\widehat{P}(t)|^2.
\end{equation} Then, by (\ref{bfo}), (\ref{convmass1}),  Plancherel's Identity, and the fact that the Fourier transform takes convolutions to products (which in this finite context follows easily from orthogonality of characters), we have \begin{equation} \label{convmass2} \theta\sigma^2 \leq \frac{4L^2}{X^2} \sum_{t\in \Z_L }|\widehat{f_B}(t)|^2|\widehat{P}(t)|^2=\frac{4L}{X^2} \sum_{x\in \Z_L } |f_B*\widetilde{P}(x)|^2 \end{equation} where  $\widetilde{P}(n)= P(-n)$,   
\begin{equation*}\label{conv}L \left( f_B * \widetilde{P}(x)\right)=\sum_{y\in \Z_L}f_B(x)P(x-y)=|B\cap(P+x)|-\sigma X,
\end{equation*}
and $P+x=\{x+\ell q : 1\leq \ell \leq X\}$. We now take advantage of the fact that $f_B$, and consequently $f_B * \widetilde{P}$, has mean value zero. In other words,
\begin{equation} \label{mv0} \sum_{n\in\Z} f_B * \widetilde{P}(n)=0.
\end{equation} As with any real valued function, we can write \begin{equation} \label{pospart} |f_B * \widetilde{P}|= 2(f_B * \widetilde{P})_+ -f_B * \widetilde{P},\end{equation} where $(f_B * \widetilde{P})_+=\max\{f_B * \widetilde{P},0\}$. 

\noindent For a moment, let us suppose that $f_B * \widetilde{P}(x)>\sigma X/L$ for some $x \in \Z_L.$ We note that $P+x$ either lifts to a genuine arithmetic progression of integers, or the disjoint union of two such progressions. In the former case, $P+x$ satisfies the conclusion of the lemma. In the latter case, if we  call the two progressions $P_1$ and $P_2$ and their lengths $X_1$ and $X_2$, then we have $(|B\cap P_1|-\sigma X_1)+(|B\cap P_2|-\sigma X_2) \geq \sigma X$. Consequently, at least one of the two progressions has length at least $\sigma X /2 $ and satisfies the required relative density condition.  

\noindent With this observation made, we can now assume that $f_B * \widetilde{P}(x)\leq \sigma X/L$ for all $x\in \Z_L$, which combined with the trivial bound of $f_B * \widetilde{P}(x)\geq -\sigma X/L$ yields \begin{equation} \label{linfty} |f_B * \widetilde{P}(x)| \leq \sigma X/L \quad \text{for all} \quad x \in \Z_L.
\end{equation}
  
\noindent Then, by (\ref{convmass2}), (\ref{mv0}), (\ref{pospart}), and (\ref{linfty}), we have
\begin{equation}\label{posmass} \sum_{x\in\Z_L} (f_B * \widetilde{P})_+(x) = \frac{1}{2} \sum_{x\in \Z_L} |f_B * \widetilde{P}(x)| \geq \frac{L}{2\sigma X} \sum_{x \in \Z_L} |f_B * \widetilde{P}(x)|^2\geq\frac{\theta\sigma X}{8}.  
\end{equation}  

\noindent We now let $E$ denote the reduced residues $x\in\Z_L$ such that $x+Xq$, as an integer, is less than $L$. Crucially, this means that if $x\in E$, then $P+x$ does not ``wrap around" and in fact lifts to a genuine arithmetic progression of integers. We see that $|\Z_L \setminus E|\leq qX$, we have assumed $f_B * \widetilde{P}$ is bounded above by $\sigma X/L$,  and by our choice of $X$ we have $qX \leq \theta L/16$,  hence \begin{equation*} \label{EF} \sum_{x\in E} (f_B * \widetilde{P})_+(x)  \geq \frac{\theta\sigma X}{8}- \frac{\sigma q X^2}{L} \geq \frac{\theta \sigma X}{16}.  
\end{equation*}

\noindent Noting that $|E|\leq L$ and recalling that $L\left(f_B * \widetilde{P}(x)\right)= |B\cap (P+x)| - \sigma X$, we have that there exists $x\in \Z$ such that the genuine arithmetic progression $P+x \subseteq \Z$ satisfies $$|B\cap(P+x)| \geq \sigma(1 + \theta/16)X,$$ as required.
\end{proof}

Noting that the progression $P$ in the conclusion of Lemma \ref{dinczn} can be partitioned into subprogressions of step size $\lambda(q)\leq q^k$, Lemma \ref{dinczn} and Proposition \ref{inh} immediately combine to yield the following iteration lemma, corresponding to Case 1 discussed in the overview, from which we deduce Lemma \ref{outer}.

\begin{lemma} \label{inc} Suppose $B\subseteq [1,L]$ with $|B|=\sigma L$ and $(B-B) \cap I(h_d) =\emptyset$. If  
\begin{equation}\label{Bmass} \sum_{t\in \mathbf{M}_q(L,\cQ)} |\widehat{B}(t)|^2 \geq \sigma^2(\log N)^{-1+\epsilon},\end{equation}
 for some $q \leq \cQ$, then there exists $B'\subseteq [1,L']$ satisfying $L' \gg \sigma L/\cQ^{k+1},$ $$(B'-B')\cap  I\left(h_{qd}\right)=\emptyset,$$ and
\begin{equation*} 
|B'|/L'\geq  \sigma(1+(\log N)^{-1+\epsilon}/16).\end{equation*}

\end{lemma} 

\subsection*{Proof of Lemma \ref{outer}}

Setting $A_0$=$A$, $N_0=N$, $\delta_0=\delta$, and $d_0=1$, we iteratively apply Lemma \ref{inc}. This yields, for each $j$, a set $A_j \subseteq [1,N_j]$ with $|A_j|=\delta_jN_j$ and $$(A_j-A_j)\cap I(h_{d_j})=\emptyset,$$satisfying
\begin{equation} \label{Ndeld} N_j \geq \delta^j N/(C\cQ)^{(k+1)j}, \quad \delta_j \geq \delta_{j-1}(1+(\log N)^{-1+\epsilon}/16),  \quad d_j \leq \cQ^{j},
\end{equation}
where $C$ is an absolute constant, as long as either
\begin{equation}\label{newmass}\max_{q\leq \cQ}\sum_{t\in \mathbf{M}_q(L,\cQ)} |\widehat{A_j}(t)|^2\geq \delta_j^2(\log N)^{-1+\epsilon}
\end{equation}
or $$|A_j\cap[1,N_j/2]|<\delta_jN_j/3,$$ 

\noindent as the latter condition implies $A_j$ has density at least $3\delta_j/2$ on the interval $(N_j/2,N_j]$, in which case we can take $A_{j+1}$ to be a shift of $A_j\cap (N_j/2,N_j]$.

\noindent We see that by (\ref{behrend}) and (\ref{Ndeld}), the density $\delta_j$ will exceed 1 after \begin{equation*}32\log(\delta^{-1})(\log N)^{1-\epsilon} \leq (\log N)^{1-\epsilon/2}\end{equation*} steps, hence  (\ref{newmass}) fails and $|A_j\cap[1,N_j/2]|\geq \delta_jN_j/3$ for some 
\begin{equation} \label{kbound} j\leq (\log N)^{1-\epsilon/2}.\end{equation}
However, we see that (\ref{behrend}), (\ref{Ndeld}), and (\ref{kbound}) imply 
 \begin{equation*}N_j \geq \delta^{(\log N)^{1-\epsilon/2}} N/(C\cQ)^{(k+1)(\log N)^{1-\epsilon/2}} \geq Ne^{-(\log N)^{1-\epsilon/8}} \geq N^{1-\rho}, \end{equation*}
so we set $B=A_j$, $L=N_j$, $\sigma = \delta_j$, and $d=d_j$, and we see further that
\begin{equation*} d \leq \cQ^{(\log N)^{1-\epsilon/2}}\leq e^{(\log N)^{1-\epsilon/4}} \leq N^{\rho/2}, \end{equation*} 
as required.
\qed

\

Our task for Section \ref{mainsec} is now completely reduced to a proof of Lemma \ref{itn}. 
\subsection{The Inner Iteration: Proof of Lemma \ref{itn}}\label{itn2} Let $B\subseteq [1,L]$ and $d\in \N$ be as in the conclusion of Lemma \ref{outer}, let $B_1=B\cap[1,L/2]$, and suppose $\sigma =|B|/L \geq \cQ^{-1/m}$. Further, let $M=\lfloor (L/3b_d)^{1/k}\rfloor$, where $b_d$ is the leading coefficient of $h_d$. Letting $$H=\{n \in \N : 0<h_d(n)<L/3\}$$  we note that by (\ref{symBIG}) we have \begin{equation}\label{symdif2} |H \ \triangle  \ [1,M]| \ll_{h} 1.  \end{equation}  Suppose we have a set $P$ with parameters $U,V,K$ as specified in the hypotheses of Lemma \ref{itn}. We fix an element $s \in P$, we let $\eta=c_0\sigma/U$ for a sufficiently small constant $c_0=c_0(h)>0$, and we let $Y=\eta^{-(k+\epsilon)}$.

\noindent Since $(B-B) \cap I(h_d)=\emptyset$, we see that there are no solutions to \begin{equation*}a-b\equiv h_d(n) \mod{L}, \quad a\in B, \  b \in B_1, \  n\in H. \end{equation*}  Therefore, the only contributions to the sum $$\frac{1}{wL^{2}}\sum_{\substack{x\in\Z_L \\ n \in [1,M]\cap W_d\left(Y\right)}} h_d'(n)B(x+h_d(n))B_1(x)e^{2\pi \i xs/L}$$ come when $n$ lies in $H \triangle [1,M]$, which we know from (\ref{symdif2}) has size $O_h(1)$. 
Bounding the contributions from $H \triangle [1,M]$ trivially and combining this observation with the orthogonality relation $$\frac{1}{L}\sum_{t\in\Z_L} e^{2\pi \i xt/L}=\begin{cases}1 &\text{if }x=0 \\ 0 &\text {if }x\in \Z_L \setminus \{0\} \end{cases},$$ we have
\begin{equation*} \sum_{t \in \Z_L} \widehat{B}(t)\bar{\widehat{B_1}(s+t)}S(t) =\frac{1}{wL^{2}}\sum_{\substack{x\in\Z_L \\ n \in [1,M]\cap W_d\left(Y\right)}} h_d'(n)B(x+h_d(n))B_1(x)e^{2\pi \i xs/L}= O_h((wM)^{-1}),\end{equation*} 
where $$w=\prod_{p\leq Y}\left(1-\frac{j_d(p)}{p^{\gamma_d(p)}}\right) $$ and $$S(t)=\frac{1}{wL}\sum_{\substack{n=1 \\ n\in W_d\left(Y\right)}}^{M} h_d'(n)e^{2\pi \i h_d(n) t/L},$$  which immediately yields  
\begin{align*}\sum_{t\in \Z_{L} \setminus \{0\}} |\widehat{B}(t)||\widehat{B_1}(s+t)||S(t)| \geq \left|\sum_{t\in \Z_L \setminus \{0\}} \widehat{B}(t)\bar{\widehat{B_1}(s+t)}S(t)\right| = \widehat{B}(0)\left|\widehat{B_1}(s)\right|S(0) - O_h\left((wM)^{-1}\right).\end{align*} 
Therefore, since $\widehat{B}(0)=\sigma$, $|\widehat{B_1}(s)| \geq \sigma/U$, $S(0)\geq 1/4$, and $\sigma^{-1}, U\leq \cQ$, we have that  
\begin{equation}\label{massW}\sum_{t\in \Z_L \setminus \{0\}} |\widehat{B}(t)||\widehat{B_1}(s+t)||S(t)| \geq \frac{\sigma^2}{5U},
\end{equation} where we also use (\ref{logk}) to absorb the error term. 

\noindent Using a variety of exponential sum estimates, both old and new, we find that 
\begin{equation}\label{Wmin} |S(t)| \leq  \frac{\sigma}{10U} \quad \text{for all } t \in \mathfrak{m}\left(L,\eta^{-1},\eta^{-(2+\epsilon)}\right), \end{equation} 
provided we choose $c_0$ sufficiently small, and \begin{equation} \label{Wmaj} |S(t)| \ll_{h}C^{\omega(q)} q^{-1/2}\min\left\{1, \left(L\left|t/L-a/q\right|\right)^{-1} \right\} \end{equation} 
if $t \in \mathbf{M}_{a/q}(L,\eta^{-1}),$ $(a,q)=1$, and $q\leq \eta^{-(2+\epsilon)}$, where $C=C(k)$ and  $\omega(q)=|\{p \text{ prime}: p \mid q\}|$. 

\noindent \textit{Remark.} As discussed in Section \ref{sievesecNew}, the inability to establish the exponent of $-1/2$ in (\ref{Wmaj}), which we obtain thanks to our careful sieving of inputs and the application of Hensel's Lemma, had previously been the fundamental obstacle in extending this method. We highlight the need for this exponent at the end of the proof and discuss these estimates in detail in Section \ref{expest}.

\noindent We have by (\ref{Wmin}), Cauchy-Schwarz, and Plancherel's Identity that 
\begin{equation*} \sum_{t \in \mathfrak{m}\left(L,\eta^{-1},\eta^{-(2+\epsilon)}\right) } |\widehat{B}(t)||\widehat{B_1}(t)||S(t)|  \leq \frac{\sigma^2}{10U},  \end{equation*} 
which together with (\ref{massW}) yields 
\begin{equation} \label{Mmass} \sum_{t \in \mathfrak{M}\left(L,\eta^{-1},\eta^{-(2+\epsilon)}\right)  } |\widehat{B}(t)||\widehat{B_1}(t)||S(t)|  \geq \frac{\sigma^2}{10U}.  \end{equation} 
We now wish to assert that we can ignore those frequencies in the major arcs at which the transform of $B$ or $B_1$ is particularly small. In order to make this precise, we first need to invoke a sieved, weighted version of known estimates on the higher moments of Weyl sums. 

\noindent Specifically, we have that  
\begin{equation} \label{moment} \sum_{t \in \Z_L} |S(t)|^{m} \leq  C, \end{equation} 
where $C=C(m)$.

\noindent \textit{Remark on (\ref{moment}) and the choice of $m>k^2+k$.} Without any sieving of inputs, Theorem 1.1 of \cite{bdg} and the proof of Proposition 3.3 in \cite{LM2} give the desired high moment estimate for $m=k^2+k$, with a small power loss. As shown in the ``$\epsilon$ removal argument" in Section 5 of \cite{bdg}, one can apply major and minor arc estimates and eliminate the small power loss by raising the exponent by any amount. The only additional insight required here is that this same trick can be used, with our major and minor arc estimates established in Section \ref{expest}, to also eliminate any potential loss caused by the sieving of inputs.

\noindent Choosing a constant $0<c_1<(40C^{1/m})^{-m/2}$, where $C$ comes from (\ref{moment}), we define 
\begin{equation} \label{XY} \mathcal{X} = \left\{ t\in \mathfrak{M}\left(L,\eta^{-1},\eta^{-(2+\epsilon)}\right)  :\min\left\{|\widehat{B}(t)|,|\widehat{B_1}(s+t)|\right\} \leq c_1\sigma^{(m+1)/2}U^{-m/2}\right\}\end{equation}
and \begin{equation*}\mathcal{Y}=\mathfrak{M}\left(L,\eta^{-1},\eta^{-(2+\epsilon)}\right)\setminus \mathcal{X}. \end{equation*} 
Using H\"{o}lder's Inequality to exploit the higher moment estimate on $S$, followed by Plancherel's Identity, we see that 
\begin{align*}\sum_{t \in \mathcal{X}}|\widehat{B}(t)||\widehat{B_1}(s+t)||S(t)| &\leq \left(\sum_{t \in \mathcal{X}} |\widehat{B}(t)|^{\frac{m}{m-1}}|\widehat{B_1}(s+t)|^{\frac{m}{m-1}} \right)^{\frac{m-1}{m}} \left( \sum_{t \in \Z_L} |S(t)|^m \right)^{\frac{1}{m}}\\ &\leq \frac{c_1^{2/m}\sigma^{\frac{m+1}{m}}}{U}\left(\sum_{t \in \Z_L} \max\left\{|\widehat{B}(t)|^2,|\widehat{B_1}(s+t)|^2\right\}\right)^{\frac{m-1}{m}}\cdot C^{1/m} \\ & \leq \frac{\sigma^{\frac{m+1}{m}}}{40U}\left(\sum_{t \in \Z_L} |\widehat{B}(t)|^2+|\widehat{B_1}(s+t)|^2\right)^{\frac{m-1}{m}} \\ & \leq \frac{\sigma^2}{20U},\end{align*}
and hence by (\ref{Mmass}) we have 
\begin{equation}\label{Ymass}  \sum_{t \in \mathcal{Y} } |\widehat{B}(t)||\widehat{B_1}(s+t)||S(t)|  \geq \frac{\sigma^2}{20U}. \end{equation} 
For $i,j,\ell \in \N$, we define 
\begin{equation*} \mathcal{R}_{i,j,\ell} = \left\{ a/q: (a,q)=1, \ 2^{i-1} \leq q \leq 2^i,  \text{ } \frac{\sigma}{2^j} \leq \max |\widehat{B}(t)| \leq \frac{\sigma}{2^{j-1}},\text{ } \frac{\sigma}{2^{\ell}} \leq \max |\widehat{B_1}(s+t)| \leq \frac{\sigma}{2^{\ell-1}} \right\}, \end{equation*}
where the maximums are taken over nonzero frequencies  $t\in \mathbf{M}_{a/q}(L,\eta^{-1})$. We see that we have 
\begin{equation}\label{rij}\sum_{a/q\in \mathcal{R}_{i,j,\ell}}\sum_{t\in \mathbf{M}_{a/q}(L,\eta^{-1})\setminus\{0\}} |\widehat{B}(t)||\widehat {B_1}(s+t)||S(t)| \ll |\mathcal{R}_{i,j,\ell}|\frac{\sigma^2}{2^{j}2^k}\max_{a/q \in \mathcal{R}_{i,j,\ell}}\sum_{t\in \mathbf{M}_{a/q}(L,\eta^{-1})} |S(t)|. \end{equation}
Examining (\ref{Wmaj}) more closely, we see that the bound of $1$ inside the minimum will be used for at most two values of $t$, and for the remaining values of $t$ that minimum can be bounded sequentially by $1/2, 1/3, \dots, \eta$. In particular, if we sum over all $t\in \mathbf{M}_{a/q}(L,\eta^{-1})$, that minimum contributes a total of $\ll \log(\eta^{-1})$. 

\noindent Therefore, it follows from (\ref{Wmaj}), the bound $U,\sigma^{-1}\leq \cQ^{1/m}$, and the standard estimate $\omega(q)\ll \log q/\log\log q,$ that if $(a,q)=1$ and $q\leq \eta^{-(2+\epsilon)}$, then
\begin{align*}\sum_{t\in \mathbf{M}_{a/q}(L,\eta^{-1})} |S(t)| &\ll_{h} C^{\omega(q)}q^{-1/2}\log(\cQ) \\ & \ll_{h} q^{-1/2}(\log N)^{\epsilon}.\end{align*} 

\noindent Therefore, by (\ref{rij})  we have 
\begin{equation} \label{rijs}\sum_{a/q\in \mathcal{R}_{i,j,\ell}}\sum_{t\in \mathbf{M}_{a/q}(L,\eta^{-1})\setminus\{0\}} |\widehat{B}(t)||\widehat{B_1}(s+t)||S(t)| \ll_{h} |\mathcal{R}_{i,j,\ell}|\frac{\sigma^2}{2^{j}2^{\ell}}2^{-i/2}(\log N)^{\epsilon}. \end{equation}   
 
\noindent By our definitions, the sets $\mathcal{R}_{i,j,\ell}$ exhaust $\mathcal{Y}$ by taking $1\leq 2^i \leq \eta^{-(2+\epsilon)}$ and $1 \leq 2^j,2^{\ell} \leq U^{m/2}/c_1\sigma^{(m-1)/2}$, a total search space of size $\ll (\log \cQ)^3 $. Therefore, by (\ref{Ymass}) and (\ref{rijs}) there exist $i,j,\ell$ in the above range with 
\begin{equation*}\frac{\sigma^2}{U(\log \cQ)^3} \ll_{h}|\mathcal{R}_{i,j,\ell}|\frac{\sigma^2}{2^j2^{\ell}}2^{-i/2}(\log N)^{\epsilon}. \end{equation*}
In other words, we can set $V_s=2^{i}$, $W_s=2^j$, and $U_s=2^{\ell}$ and take an appropriate nonzero frequency from each of the pairwise disjoint major arcs specified by $\mathcal{R}_{i,j,\ell}$ to form a set 
\begin{equation*} P_s \subseteq \left\{t \in \bigcup_{q=V_s/2}^{V_s} \mathbf{M}_{q}(L,\eta^{-1}):  \text{ } |\widehat{B_1}(s+t)| \geq \frac{\sigma}{U_s} \right\} \end{equation*} 
which satisfies 
\begin{equation}\label{Ps} |P_s| \gg_{h} \frac{U_sW_sV_s^{1/2}}{U(\log N)^{2\epsilon}},  \quad |P_s\cap\mathbf{M}_{a,q}(L,\eta^{-1})|\leq 1 \text{ \ whenever \ } q \leq V_s,  \end{equation} and 
\begin{equation}\label{maxBW} \max_{ t \in \mathbf{M}_{a/q}(L,\eta^{-1})\setminus \{0\}} |\widehat{B}(t)| \geq \frac{\sigma}{W_s} \text{ whenever } q\leq V_s \text{ and }  \mathbf{M}_{a/q}(L,\eta^{-1})\cap P_s \neq \emptyset,
\end{equation} 

\noindent noting by disjointness that $a/q \in \mathcal{R}_{i,j,\ell}$  whenever $q\leq V_s$ and  $\mathbf{M}_{a/q}(L,\eta^{-1})\cap P_s \neq \emptyset $. 

\noindent We now observe that by pigeonholing there is a subset $\tilde{P} \subseteq P$ with $|\tilde{P}| \gg |P|/(\log \cQ)^3$, and hence
\begin{equation}\label{tilde}|\tilde{P}| \gg |P|/(\log N)^{\epsilon},\end{equation}  
for which the triple $U_s,W_s,V_s$ is the same. We call  those common parameters $\tilde{U},\tilde{W}$ and $\tilde{V}$, respectively, and we can now foreshadow by asserting that the claimed parameters in the conclusion of Lemma \ref{itn} will be $U'=\tilde{U}$, $V'=\tilde{V}V$, and $K'=K+\eta^{-1}$, which do satisfy the purported bound. 

\noindent We let 
\begin{equation*} \mathcal{R} = \left\{\frac{a}{q}+\frac{b}{r}:  s\in \mathbf{M}_{a/q}(L,K) \text{ for some } s \in \tilde{P} \text{ and }  t \in \mathbf{M}_{b/r}(L,\eta^{-1}) \text{ for some } t \in P_s \right\}. \end{equation*} 

\noindent By taking one frequency $s+t$ associated to each element in $\mathcal{R}$,  we form our set $P'$, which immediately satisfies conditions (\ref{P'1}) and (\ref{P'2}) from the conclusion of Lemma \ref{itn}. However, the crucial condition (\ref{P'3}) on $|P'|$, which by construction is equal to $|\mathcal{R}|$, remains to be shown. To this end, we invoke the work on the combinatorics of rational numbers found in \cite{PSS} and \cite{BPPS}.

\begin{lemma}[Lemma CR of \cite{BPPS}]\label{CR} \begin{equation*} |\mathcal{R}|\geq \frac{|\tilde{P}|(\min_{s\in \tilde{P}}|P_s|)^2}{\tilde{V}E\tau^8(1+\log V)},  \end{equation*} where \begin{equation*} E= \max_{r\leq \tilde{V}} \Bigl|\Bigl\{ b  : \ (b,r)=1, \  \mathbf{M}_{b/r}(L,\eta^{-1})\cap \bigcup_{s\in \tilde{P}} P_s\neq \emptyset \Bigr\}\Bigr|, \end{equation*} $\tau(q)$ is the divisor function and $\tau = \max_{q\leq V\tilde{V}} \tau(q)$.
\end{lemma}

\noindent It is a well-known fact of the divisor function that $\tau(n) \leq n^{1/\log\log n}$ for large $n$, and since $\eta^{-1}, V\tilde{V} \leq \cQ$, we have that $\tau \leq (\log N)^{\epsilon}$. We also have from (\ref{maxmass}) that 
\begin{equation} \label{1mnec} \sigma^2(\log N)^{-1+\epsilon} \geq \max_{r \leq \cQ} \sum_{t \in \mathbf{M}_r(L,\cQ)} |\widehat{B}(t)|^2 \geq \max_{r \leq \tilde{V}} \sum_{t \in \mathbf{M}_r(L,\eta^{-1})} |\widehat{B}(t)|^2 \geq \frac{\sigma^2}{\tilde{W}^2}E,\end{equation}
where the last inequality follows from (\ref{maxBW}),  and hence $E \leq \tilde{W}^2(\log N)^{-1 + \epsilon}.$

 \noindent Combining the estimates on $\tau$ and $E$ with (\ref{Ps}), (\ref{tilde}), and Lemma \ref{CR}, we have\begin{equation} \label{sqrtnec} |P'| \gg_{h} \frac{|P|}{(\log N)^{\epsilon}} \frac{\tilde{U}^2\tilde{W}^2\tilde{V}}{U^2(\log N)^{4\epsilon}} \frac{(\log N)^{1-\epsilon}}{\tilde{V}\tilde{W}^2(\log N)^{9\epsilon}} = \tilde{U}^2\frac{|P|}{U^2}(\log N)^{1-15\epsilon}.\end{equation} \\
Recalling that we set $U'=\tilde{U}$, we see that for sufficiently large $N$ we have $$\frac{|P'|}{(U')^2} \geq \frac{|P|}{U^2}(\log N)^{1-16\epsilon},$$ as claimed. \qed

\noindent \textit{Remark on the necessity of square root cancellation.} We now retrospectively observe that the exponent $-1/2$ in $(\ref{Wmaj})$ directly leads to the exponent $1/2$ on $\tilde{V}$ in (\ref{Ps}). That factor of $\tilde{V}^{1/2}$ is then squared via Lemma \ref{CR}, which then cancels the factor of $\tilde{V}$ in the denominator in (\ref{sqrtnec}). As the logarithmic gain in $|P'|/(U')^2$ is the crux of the method, the argument completely breaks down without this full cancellation of $\tilde{V}$. 

\noindent \textit{Remark on the hypothesis $\max\{U,V,K,\sigma^{-1}\} \leq Q^{1/m}$ in Lemma \ref{itn}.} We make this assumption for convenience, as it ensures that not only the original parameters $U,V,K$, but also the new parameters $U',V',K'$, are bounded above by $\cQ$. This precise upper bound is only strictly crucial for the intermediate parameter $\tilde{V}$, as it is required to apply (\ref{maxmass}) in (\ref{1mnec}). This hypothesis can be carefully softened, though doing so does not lead to a stronger final result.

\noindent \textit{Remark on possible values of $c$ in Theorem \ref{main1}.} In examining the proof of Theorem \ref{main1}, we see that the value of $c$ in that theorem is determined by the exponent $m/2$ on $\max\{U,V,K\}$ in the conclusion of Lemma \ref{itn}. In examining the proof of Lemma \ref{itn}, we see that this particular exponent is necessitated by the application of the high moment estimate (\ref{moment}), which leads to the definition of $\mathcal{X}$ in (\ref{XY}). The quoted result from \cite{bdg} makes our choice of $m>k^2+k$ sharp, so for our chosen proof technique, our range for $c$ is optimal. However, in \cite{BPPS}, which restricts to $h(n)=n^k$, the constant $c=1/4$ is obtained independent of $k$, and with slight modification can be made arbitrarily close to $1/\log3$, as is done in \cite{HLR}. However, this seems to be a hard limit of the method, as even if the exponent on $\max\{U,V,K\}$ in Lemma \ref{itn} could be taken independent of the high moment estimate, we would still be limited by the upper bound $V' \ll V\eta^{-(2+\epsilon)} \ll VU^{2+\epsilon}\sigma^{-(2+\epsilon)}$ seen in the proof, which requires that exponent to be at least $3$.

\section{Exponential Sum Estimates}\label{expest} In this section, we either invoke or prove all exponential sum estimates necessary to establish the crucial major and minor arc upper bounds in Section \ref{itn2}, namely (\ref{Wmin}), and (\ref{Wmaj}).  Throughout this section, given a polynomial $g\in \Z[x]$, we let  $\gamma(p)$, $j(p)$, and $W(Y)$ be defined in terms of $g$ as in Section \ref{sievesec}. The first lemma provides asymptotic formulae for the relevant sifted Weyl sums near rationals with small denominator. 

\begin{lemma}\label{SasymW} Suppose $k\in \N$, $g(x)=a_0+a_1x+\cdots+a_kx^k \in \Z[x]$, and let $J=|a_0|+|a_1|+\cdots+|a_k|$. If $X,Y > 0$, $a,q\in \N$, $\alpha=a/q+\beta$, and $c\log (X/q) \geq \log Y \log \log Y$, then \begin{align*}\sum_{\substack{n=1 \\ n \in W(Y)}}^Xg'(n)e^{2\pi \textnormal{i}g(n)\alpha}&=\frac{1}{q} \prod_{\substack{ p\leq Y \\ p^{\gamma(p)}\nmid q}}\left(1-\frac{j(p)}{p^{\gamma(p)}}
\right)\sum_{\substack{s=0 \\ s \in W^{q}(Y)}}^{q-1}e^{2\pi \textnormal{i} g(s)\frac{a}{q}}\int_0^Xg'(x)e^{2\pi\textnormal{i}g(x)\beta}dx\\\\& + O\left(kJX^ke^{-c\frac{\log\left(\frac{X}{q}\right)}{\log Y}}(1+JX^k|\beta|)\right),\end{align*} where $c=c(k,\textnormal{cont(g)})>0$.

\end{lemma}
 
\begin{proof} We begin by noting that for any $a,q \in \N$ and $0\leq x \leq X$,

\begin{align*}\label{rat2} \sum_{\substack{n=1 \\ n \in W(Y)}}^x g'(n)e^{2\pi \i g(n)a/q}&=\sum_{s=0}^{q-1}e^{2\pi \i g(s)a/q}\sum_{\substack{n=1 \\ n\in W(Y) \\n\equiv s \text{ mod }q}}^x g'(n) \\\\ &= \frac{g(x)}{q} \prod_{\substack{ p\leq Y \\ p^{\gamma(p)}\nmid q}}\left(1-\frac{j(p)}{p^{\gamma(p)}}\right) \sum_{\substack{s=0 \\ s \in W^{q}(Y)}}^{q-1}e^{2\pi \i g(s)a/q}+ O\left(kJX^ke^{-c\frac{\log (X/q)}{\log Y}} \right),\end{align*}
\

\noindent since for $s\in W^q(Y)$ we have by the same calculation as Proposition \ref{brunprop} and partial summation that \begin{equation*} \label{count} \sum_{\substack{n=1 \\ n\in W(Y) \\n\equiv s \text{ mod }q}}^x g'(n)=\frac{g(x)}{q}\prod_{\substack{ p\leq Y \\ p^{\gamma(p)}\nmid q}}\left(1-\frac{j(p)}{p^{\gamma(p)}}\right) + O\left(\frac{kJX^k}{q}e^{-c\frac{\log (X/q)}{\log Y}} \right),\end{equation*} whereas for $s\notin W^q(Y)$ the sum is zero. We note that $c>0$ depends only on $k$ and cont$(g)$ by Prop. \ref{idzero}.

\noindent Then, by partial summation we have that if $\alpha=a/q+\beta$, then
\begin{align*}\sum_{\substack{n=1 \\ n \in W(Y)}}^Xg'(n)e^{2\pi \textnormal{i}g(n)\alpha} &= \frac{1}{q}\prod_{\substack{ p\leq Y \\ p^{\gamma(p)}\nmid q}}\left(1-\frac{j(p)}{p^{\gamma(p)}}\right)\sum_{\substack{s=0 \\ s \in W^{q}(Y)}}^{q-1}e^{2\pi \i g(s)a/q}\\\\& \cdot \left(g(X)e^{2\pi \i g(X)\beta}-\int_0^Xg(x)(2\pi \i \beta g'(x))e^{2\pi \i g(x)\beta}dx\right)\\\\ & +O\left(kJX^ke^{-c\frac{\log\left(\frac{X}{q}\right)}{\log Y}}(1+JX^k|\beta|)\right).
\end{align*} 

\noindent Finally, noting that  $$g(X)e^{2\pi \i g(X)\beta}-\int_0^Xg(x)(2\pi \i \beta g'(x))e^{2\pi \i g(x)\beta}dx= \int_0^Xg'(x)e^{2\pi\textnormal{i}g(x)\beta}dx, $$ the lemma follows.
\end{proof} 

To establish the required square root cancellation for the restricted exponential sums that arise in the conclusion of Lemma \ref{SasymW}, we use the following standard fact about lifting roots of polynomials, the details of which motivate the sieve outlined in Section \ref{sievesec}. Here we focus only on the existence of the lifts, ignoring the extent of uniqueness or ``closeness", and this statement follows, for example, from Proposition 2 in Chapter II, Section 2 of \cite{Lang}.
 
\begin{lemma}[Hensel's Lemma] \label{hensel} Suppose $g\in \Z[x]$, $p$ is prime, and $n,\gamma,j\in \N$ with $j\geq 2\gamma-1$. If $$g(n)\equiv 0 \ \textnormal{mod} \ p^{2\gamma-1}$$ and $g'(n)\not\equiv 0$ \textnormal{mod} $p^{\gamma}$, then there exists $m\in \N$ with $g(m) \equiv 0$ \textnormal{mod} $p^j$.
\end{lemma}

Armed with Lemma \ref{hensel}, we exhibit the restricted exponential sum estimate that was previously the missing ingredient to proving Theorem \ref{main1}.

\begin{lemma}\label{exponentialest} If $g\in \Z[x]$ with $\deg(g)=k\geq 2$, $a,q\in \N$ with $(a,q)=1$, and $Y>0$, then \begin{equation*}\left|\sum_{\substack{s=0 \\ s \in W^{q}(Y)}}^{q-1}e^{2\pi \textnormal{i} g(s)a/q}\right| \ll_k \gcd(\textnormal{cont}(g),q)^{3} C^{\omega(q)}\begin{cases} q^{1/2} &\text{if }q\leq Y \\ q^{1-1/k} &\text{for all }q \end{cases}, \end{equation*}where $C=C(k)$ and $\omega(q)$ is the number of distinct prime factors of $q$.
\end{lemma}
\begin{proof} Factor $q=q_1\cdots q_4$, where $q_1,\dots,q_4$ are pairwise coprime, $q_1$ houses the prime power factors $p^j$ of $q$ satisfying $p\leq Y$, $\gamma(p)>1$, and $j<2\gamma(p)$, $q_2$ is a product of distinct primes $p\leq Y$ satisfying $\gamma(p)=1$, $q_3$ is the product of $p^j$ satisfying $p\leq Y$, $j\geq 2\gamma(p)$, and all prime factors of $q_4$ are greater than $Y$.  

\noindent By the Chinese Remainder Theorem, we have \begin{equation*} \sum_{\substack{s=0 \\ s \in W^{q}(Y)}}^{q-1}e^{2\pi \textnormal{i} g(s)a/q}= \prod_{i=1}^{4} \sum_{\substack{s=0 \\ s \in W^{q_i}(Y)}}^{q_i-1}e^{2\pi \textnormal{i} g(s)a_i/q_i}, \end{equation*} where $a_1,\dots,a_4$ are the unique residues satisfying $$\frac{a}{q} \equiv \frac{a_1}{q_1}+\cdots+\frac{a_{4}}{q_4} \ \text{mod }1. $$
By definition of $\gamma$, $g'$ is identically zero modulo $p^{\gamma(p)-1}$, and since any common factor of the coefficients of $g'$ divides $k!\text{cont}(g)$, we have by Proposition \ref{idzero} that $p^{\gamma(p)-1} \ll_k \gcd(\text{cont}(g),p^{\gamma(p)-1}),$. Therefore, since $p^{2\gamma(p)-1}\leq p^{3(\gamma(p)-1)}$ if $\gamma(p)>1$, we have $$q_1\ll_k \gcd(\text{cont}(g),q_1)^3.$$ Further decomposing $q_2$ into a product of primes, using the fact that $g'$ has at most $k-1$ roots modulo each of these primes,  and applying the standard Weil bound (see for example Theorem 3.1 of \cite{Kow}), we have  $$\left|\sum_{\substack{s=0 \\ s \in W^{p}(Y)}}^{p-1}e^{2\pi \textnormal{i} g(s)b/p}\right| \ll_k p^{1/2} $$ provided $p\nmid b\text{cont}(g)$, and hence \begin{equation*}\label{q2} \left|\sum_{\substack{s=0 \\ s \in W^{q_2}(Y)}}^{q_2-1}e^{2\pi \textnormal{i} g(s)a_2/q_2}\right| \leq C^{\omega(q_2)}\gcd(\text{cont}(g),q_2)^{1/2} q_2^{1/2}, \end{equation*} where $C=C(k)$. Now suppose that $p$ is prime, $j\geq 2\gamma(p)$, and $\ell=2\gamma(p)-1$. If $0\leq s \leq p^j-1$ and $s_1$ is the reduced residue class of $s$ modulo $p^{\ell}$, then $g(s)\equiv p^{\ell}s_2+g(s_1)$ mod $p^j$ for some $0\leq s_2\leq p^{j-\ell}-1$. Conversely,  if $0\leq s_1 \leq p^{\ell}-1$ with  $g'(s_1)\not\equiv 0$ mod $p^{\gamma(p)}$, then for every $0\leq s_2 \leq p^{j-\ell}-1$, Lemma \ref{hensel} applied to the polynomial $g(x)-p^{\ell}s_2+g(s_1)$ yields $0\leq s \leq p^{j}-1$ with $g(s)\equiv p^{\ell}s_2+g(s_1)$ mod $p^j$.

\noindent In other words, the map $F$ on $\Z/p^{j-\ell}\Z$ defined by $g(p^{\ell}s_2+s_1)\equiv p^{\ell}F(s_2)+g(s_1)$ mod $p^j$ is a bijection.  In particular, if $p\nmid b$, then \begin{align*}\sum_{\substack{s=0 \\ s \in W^{p^j}(Y)}}^{p^j-1}e^{2\pi \textnormal{i} g(s)b/p^j} &= \sum^{p^{\ell}-1}_{\substack{s_1=0 \\ g'(s_1)\not\equiv 0 \text{ mod }p^{\gamma(p)}}} \sum_{s_2=0}^{p^{j-\ell}-1}e^{2\pi \i g(p^{\ell}s_2+s_1)b/p^j}\\\\  &=\sum^{p^{\ell}-1}_{\substack{s_1=0 \\ g'(s_1)\not\equiv 0 \text{ mod }p^{\gamma(p)}}} \sum_{s_2=0}^{p^{j-\ell}-1}e^{2\pi \i \left(p^{\ell}s_2+g(s_1)\right)b/p^j} \\\\ &=0,\end{align*} where the last equality is the fact that the sum in $s_2$ runs over the full collection of $p^{j-\ell}$-th roots of unity. Therefore, we have that \begin{equation*}\sum_{\substack{s=0 \\ s \in W^{q_3}(Y)}}^{q_3-1}e^{2\pi \textnormal{i} g(s)a_3/q_3}=\begin{cases}1 & \text{if }q_3=1 \\ 0 &\text{else} \end{cases}. \end{equation*}  Finally, noting that $W^{q_4}(Y)=\N$, we utilize the  standard complete sum estimate (see \cite{Chen} for example)  \begin{equation*}\left|\sum_{s=0}^{q_4-1}e^{2\pi \textnormal{i} g(s)a_4/q_4}\right| \ll_k \gcd(\text{cont}(g),q_4)^{1/k}q_4^{1-1/k}, \end{equation*} and the lemma follows.
\end{proof} 

\

We now invoke a variation of the most traditional minor arc estimate, Weyl's Inequality.

\

\begin{lemma}[Lemma 3 in \cite{CLR}] \label{weyl2}  Suppose $k\in \N$, $g(x)=a_0+a_1x+\cdots+a_{k}x^{k}$ with $a_1\dots,a_k \in \R$ and $a_{k} \in \N$. If $X>0$, $a,q\in \N$ with $(a,q)=1$, and $|\alpha-a/q|<q^{-2}$, then 
$$ \left|\sum_{n=1}^X e^{2\pi \i g(n)\alpha} \right| \ll_{k} X \left(a_{k}\log^{k^2}(a_{k}qX)\left(q^{-1}+X^{-1}+\frac{q}{a_kX^k}\right) \right)^{2^{-k}}.$$
\end{lemma}

\

Finally, we carefully adapt Lemma \ref{weyl2} to our particular sieve to get the desired estimates far from rationals with  small denominator.

\

\begin{lemma}\label{weyl3} Suppose $k\in \N$, $g(x)=a_0+a_1x+\cdots+a_kx^k \in \Z[x]$  with $a_k>0$. Suppose further that $X,Y,Z\geq 2$, $YZ\leq X$, and $a,q\in \N$ with $(a,q)=1$. If $|\alpha-a/q|<q^{-2}$, then $$\left|\sum_{\substack{n=1 \\ n \in W(Y)}}^X e^{2\pi \i g(n)\alpha} \right| \ll_{k} \textnormal{cont}(g)^6(\log Y)^{ek} X\left(e^{-\frac{\log Z}{\log Y}}+\left(a_{k}\log^{k^2}(a_{k}qX)\left(q^{-1}+\frac{Z}{X}+\frac{qZ^k}{a_kX^k}\right) \right)^{2^{-k}} \right). $$
\end{lemma}

\

\begin{proof} Let $P$ be the set of products $p_1^{\gamma(p_1)}\cdots p_s^{\gamma(p_s)}$ with $p_1<\cdots<p_s\leq Y$. By inclusion-exclusion, \begin{equation} \label{P1P2} \left|\sum_{\substack{n=1 \\ n \in W(Y)}}^X e^{2\pi \i g(n)\alpha}\right| = \left|\sum_{D \in P} (-1)^{\omega(D)} \sum_{\substack{n=1 \\ g'(n) \equiv 0 \text{ mod }D}}^X  e^{2\pi \i g(n)\alpha} \right|
\end{equation} and further by Proposition \ref{idzero}\begin{align*}\left|\sum_{\substack{D \in P \\ D\leq Z}} (-1)^{\omega(D)} \sum_{\substack{n=1 \\ g'(n) \equiv 0 \text{ mod }D}}^X  e^{2\pi \i g(n)\alpha} \right| &\ll_k \text{cont}(g)^2\sum_{D \in P} k^{\omega(D)} \max_{0\leq b \leq D}\left| \sum_{n=0}^{X/D} e^{2\pi \i g(Dn+b)\alpha}\right|. \end{align*} This inequality uses a few things. Specifically, the inner range of summation on the left side is a disjoint union of arithmetic progressions modulo $D$, where the number of progressions is the product of the number of roots of $g'$ modulo $p_i^{\gamma(p_i)}$.  We then use that if $\gamma(p)=1$, then $g'$ has fewer than $k$ roots modulo $p$, while the $\text{cont}(g)^2$ term accounts for the primes $p$ for which $\gamma(p)>1$, as  $p^{\gamma(p)}\leq p^{2(\gamma(p)-1)} \ll_k \text{cont}(g)^2$. 

\noindent Further, we see from Lemma \ref{weyl2} that \begin{align*}\sum_{\substack{D \in P \\ D\leq Z}} k^{\omega(D)} \max_{0\leq b \leq D}\left| \sum_{n=0}^{X/D} e^{2\pi \i g(Dn+b)\alpha}\right| &\ll_k \sum_{\substack{D \in P \\ D\leq Z}} k^{\omega(D)} \frac{X}{D} \left(a_{k}\log^{k^2}(a_{k}qX)\left(q^{-1}+\frac{D}{X}+\frac{qD^k}{a_kX^k}\right) \right)^{2^{-k}} \\\\ &\ll_{k} X \left(a_{k}\log^{k^2}(a_{k}qX)\left(q^{-1}+\frac{Z}{X}+\frac{qZ^k}{a_kX^k}\right) \right)^{2^{-k}}\sum_{D\in P} \frac{k^{\omega(D)}}{D} \\\\ &\ll_k X(\log Y)^k\left(a_{k}\log^{k^2}(a_{k}qX)\left(q^{-1}+\frac{Z}{X}+\frac{qZ^k}{a_kX^k}\right) \right)^{2^{-k}},\end{align*} where the last inequality uses that if $C>0$, then \begin{equation}\label{Z}\sum_{D\in P} \frac{C^{\omega(D)}}{D} = \prod_{p\leq Y} \left(1+\frac{C}{p^{\gamma(p)}}\right) \leq   \prod_{p\leq Y} \left(1+\frac{C}{p}\right) \ll (\log Y)^C.\end{equation} For the contribution of large $D$ to (\ref{P1P2}), by considering the cancellation in the alternating term $(-1)^{\omega(D)}$, we can reduce to the case of $D\in P$ with $Z<D\leq qZ$, where $q$ is the largest prime power $p^{\gamma(p)}$ with $p\leq Y$, so in particular $q\ll_k \text{cont}(g)Y$. Then, since $D\ll_k \text{cont}(g)^2Y^{\omega(D)}$ and $Y\geq 2$, we know that  \begin{equation}\label{logQY} \text{cont}(g)^3e^{\omega(D)-\frac{\log Z}{\log Y}} \gg_k 1. \end{equation} Finally, by trivially bounding the inner sum and applying (\ref{logQY}) and (\ref{Z}), we have that if $YZ\leq X$, then \begin{align*}\left|\sum_{\substack{D \in P \\ D>Z} } (-1)^{\omega(D)} \sum_{\substack{n=1 \\ g'(n) \equiv 0 \text{ mod }D}}^X  e^{2\pi \i g(n)\alpha} \right| &  \ll_k \text{cont}(g)^2 \sum_{\substack{D \in P \\ Z<D\leq qZ}} k^{\omega(D)} \left(1+\frac{X}{D}\right) \\\\ & \ll_k \text{cont}(g)^3\sum_{\substack{D \in P \\ D> Z}}k^{\omega(D)}\frac{X}{D} \\\\ &\ll_k \text{cont}(g)^6 e^{-\frac{\log Z}{\log Y}}X \sum_{D\in P} \frac{(ek)^{\omega(D)}}{D} \\\\& \ll \text{cont}(g)^6 e^{-\frac{\log Z}{\log Y}} (\log Y)^{ek}X,\end{align*} and the estimate follows.\end{proof}

\subsection{Proof of (\ref{Wmin}) and (\ref{Wmaj})} \label{it} We return to the setting of the proof of Lemma \ref{itn} in Section \ref{itn2}, recalling all assumptions,  notation, and fixed parameters. Further, we let $$Z=e^{(\log\log N)^3},$$ noting that \begin{equation}\label{cQ2} e^{\frac{\log Z}{\log Y}}>\cQ^2. \end{equation} Fixing $t\in \Z_L$, the pigeonhole principle guarantees the existence of $1\leq q \leq L/Z^{2k}$ and $(a,q)=1$ with $$\left|\frac{t}{L}-\frac{a}{q}\right|<\frac{Z^{2k}}{qL}.$$ Letting $\beta=t/L-a/q$, we have by Lemma \ref{SasymW} and (\ref{logk}) that \begin{equation}\label{Sest} S(t)=\frac{w_q}{qL} \sum_{\substack{s=0 \\ s \in W^{q}(Y)}}^{q-1}e^{2\pi \textnormal{i} h_d(s)a/q}\int_0^Mh_d'(x)e^{2\pi\textnormal{i}h_d(x)\beta}dx+ O_h\left(e^{-c\frac{\log\left(\frac{M}{q}\right)}{\log Y}}Z^{3k}\right), \end{equation} where $$w_q=\prod_{\substack{ p\leq Y \\ p^{\gamma(p)}\mid q}}\left(1-\frac{j_d(p)}{p^{\gamma_d(p)}}\right)^{-1}. $$

\noindent We note that $w_q \ll_h 2^{\omega(q)}$ by Proposition \ref{idzero} and Lemma \ref{content}. Combining (\ref{Sest}) and Lemma \ref{exponentialest} with  \begin{equation} \label{VDC}  \left|\int_0^M h_d'(x)e^{2\pi\textnormal{i}h_d(x)\beta}dx\right|=\left|\int_0^{h_d(M)}e^{2\pi \i y \beta} \d y\right| \ll \min \{L, |\beta|^{-1} \}  \end{equation} yields  (\ref{Wmaj}) if $$q\leq \eta^{-(2+\epsilon)} \quad \text{and} \quad  |\beta|<(\eta L)^{-1},$$ as well as (\ref{Wmin}) if $$q\leq \eta^{-(2+\epsilon)} \ \text{and} \ |\beta|\geq(\eta L)^{-1} \quad  \text{or} \quad  \eta^{-(2+\epsilon)}<q\leq L^{2k\rho}.$$ 

\noindent Finally, recalling that the leading coefficient $b_d$ of $h_d$ satisfies $b_d \ll_h d^k \leq L^{k\rho}$, we have by Lemma \ref{weyl3}, Lemma \ref{content}, (\ref{cQ2}), and partial summation that if $L^{2k\rho}\leq q \leq L/Z^{2k}$, then $$|S(t)|\ll_h \cQ^{-1},$$ and in particular (\ref{Wmin}) holds. \qed

\section{Single Iteration Method: Proof of Theorem \ref{main2}} \label{singsec}

In this section, we further exploit the estimates established in Section \ref{expest} and apply a more traditional $L^2$ density increment, essentially an improved, streamlined version of S\'ark\"ozy's \cite{Sark1} original method, in order to prove Theorem \ref{main2}. The core of this method has been utilized in \cite{Lucier}, \cite{LM}, \cite{Ruz}, and \cite{Rice}, among others.

We also provide a brief discussion on sums of three or more polynomials in Section \ref{finalsec}. We begin with another preliminary discussion of the circle method, this time with a continuous frequency domain.

\subsection{Fourier analysis and the circle method on $\Z$} For this argument, rather than identify an interval of integers with a cyclic group, we embed our finite sets in $\Z$, on which we utilize an unnormalized discrete Fourier transform. Specifically, for a function $F: \Z \to \C$ with finite support, we define $\widehat{F}: \T \to \C$, where $\T$ denotes the circle  parameterized by the interval $[0,1]$ with $0$ and $1$ identified, by \begin{equation*} \widehat{F}(\alpha) = \sum_{x \in \Z} F(x)e^{-2 \pi \text{i}x\alpha}. \end{equation*}

\noindent Given $N\in \N$ and a set $A\subseteq [1,N]$ with $|A|=\delta N$,  rather than singling out the zero frequency, we examine the Fourier analytic behavior of $A$ by considering the \textit{balanced function}, $f_A$, defined by
\begin{equation*} f_A=1_A-\delta 1_{[1,N]}.\end{equation*}

\noindent We then define the major and minor arcs on $\T$, analogous to our definitions from Section \ref{ZNZ}. 

\begin{definition}Given $\gamma>0$ and $Q\geq 1$, we define, for each $q\in \N$ and $a\in [1,q]$,
$$\mathbf{M}_{a/q}(\gamma)=\left\{ \alpha \in \T : \Big|\alpha-\frac{a}{q}\Big| < \gamma \right\},$$  $$\mathbf{M}_q(\gamma)=\bigcup_{(a,q)=1} \mathbf{M}_{a/q}(\gamma),$$ and $$ \mathbf{M}'_q(\gamma)=\bigcup_{r\mid q} \mathbf{M}_r(\gamma)=\bigcup_{a=1}^q \mathbf{M}_{a/q}(\gamma).$$
We then define $\mathfrak{M}(\gamma,Q)$, the \textit{major arcs}, by
\begin{equation*} \mathfrak{M}(\gamma,Q)=\bigcup_{q=1}^{Q} \mathbf{M}_q(\gamma),
\end{equation*} 
and $\mathfrak{m}(\gamma,Q)$, the \textit{minor arcs}, by  
\begin{equation*} \mathfrak{m}(\gamma,Q)=\T\setminus \mathfrak{M}(\gamma,Q).
\end{equation*} 
We note that if $2\gamma Q^2<1$, then \begin{equation} \label{majdisj}\mathbf{M}_{a/q}(\gamma)\cap\mathbf{M}_{b/r}(\gamma)=\emptyset \end{equation}whenever  $a/q\neq b/r$ and  $q,r \leq Q$. 
\end{definition}

\subsection{Generalized Inheritance Proposition} To establish Theorems \ref{main2} and \ref{main3}, we require a generalization of Proposition \ref{inh} that replaces a single polynomial image with the sumset of a collection of polynomial images. To this end,  if $h^{(1)},\dots, h^{(\ell)}\in \Z[x]$ is a collection of intersective polynomials, then we define $\lambda_1,\dots,\lambda_{\ell}$ as in Section \ref{auxdef1} in terms of $h^{(1)},\dots,h^{(\ell)}$, respectively. Then, we define $\Lambda=\lambda_1\circ \cdots \circ \lambda_{\ell}$, where $\circ$ denotes composition, and $$\tilde{\lambda}_i=\lambda_1\circ\cdots\lambda_{i-1} \circ \lambda_{i+1}\circ \cdots \lambda_{\ell}$$ for $1\leq i \leq \ell$. We note that the $\lambda_i$ commute with each other, and in particular $\lambda_i \circ \tilde{\lambda_i}=\Lambda$. As is standard, we use $$A \pm B = \{a \pm b : a\in A, b\in B\}$$ to denote the sum and difference sets, respectively.

\begin{proposition} \label{inh2} Suppose $h^{(1)},\dots,h^{(\ell)}\in \Z[x]$ is a collection of intersective polynomials, $d_1,\dots,d_{\ell}\in \N$, and $A \subseteq \N$. If $x\in \Z$, $q\in \N$, $$(A-A)\cap \Big(I(h^{(1)}_{d_1})+\cdots+I(h^{(\ell)}_{d_{\ell}})\Big)\subseteq \{0\},$$ and $A'\subseteq \{a\in \N : x+\Lambda(q)a \in A\}$, then $$(A'-A')\cap \Big(I(h^{(1)}_{\tilde{\lambda}_1(q)d_1})+\cdots+ I(h^{(\ell)}_{\tilde{\lambda}_{\ell}(q)d_{\ell}})\Big)\subseteq \{0\}.$$ 
\end{proposition}

\begin{proof} Suppose that $A\subseteq \N$, $A'\subseteq \{a\in \N : x+\Lambda(q)a \in A\}$, and 
\begin{align*}0\neq a-a'&=\sum_{i=1}^{\ell}h^{(i)}_{\tilde{\lambda}_i(q)d_i}(n_i)=\sum_{i=1}^{\ell}\frac{h^{(i)}\left(r^{(i)}_{\tilde{\lambda}_i(q)d_i}+\tilde{\lambda}_i(q)d_in_i\right)}{\lambda_i\left(\tilde{\lambda}_i(q)d_i\right)}=\sum_{i=1}^{\ell}\frac{h^{(i)}\Big(r^{(i)}_{\tilde{\lambda}_i(q)d_i}+\tilde{\lambda}_i(q)d_in_i\Big)}{\Lambda(q)\lambda_i(d_i)}\end{align*} 
for some $n_1,\dots, n_{\ell}\in \N$, $a,a' \in A'$, with all polynomial terms having the same sign as the corresponding leading coefficient. 

\noindent By construction we know that $r^{(i)}_{\tilde{\lambda}_i(q)d_i}\equiv r^{(i)}_{d_i}$ mod $d_i$, so there exists $s_i\in \Z$ such that $r^{(i)}_{\tilde{\lambda}_i(q)d_i}=r^{(i)}_{d_i}+d_is_i$, and therefore \begin{align*}0\neq \sum_{i=1}^{\ell} h^{(i)}_{d_i}(s_i+\tilde{\lambda}_i(q)n_i)=\sum_{i=1}^{\ell}\frac{h^{(i)}(r_{d_i}^{(i)}+d_i(s_i+\tilde{\lambda}_i(q)n_i))}{\lambda_i(d_i)}=\Lambda(q)(a-a').\end{align*} Because $A'\subseteq \{a\in \N : x+\Lambda(q)a \in A\}$, we know that $\Lambda(q)(a-a')\in A-A$, hence  $$(A-A)\cap \Big(I(h^{(1)}_{d_1})+\cdots+I(h^{(\ell)}_{d_{\ell}})\Big) \not\subseteq \{0\},$$ and the contrapositive is established.
\end{proof}

\subsection{Main iteration lemma and proof of Theorem \ref{main2}} \label{mainproof} For the remainder of Section \ref{singsec} we fix intersective polynomials $g,h\in \Z[x]$, and we let $k=\deg(g)$, $\ell=\deg(h)$,  we let $D=\left(k^{-1}+\ell^{-1}\right)^{-1}$, and we let $\rho=\rho(g,h)>0$ be an appropriately small constant. Finally, for $N\in \N$ we let $$\cQ=\cQ(N)=e^{\rho(\log N)^{1/3}}.$$ 

\noindent We deduce Theorem \ref{main2}  from the following iteration lemma, which states that a set deficient in the desired pattern spawns a new, significantly denser subset of a slightly smaller interval with an inherited deficiency in the pattern associated to appropriate auxiliary polynomials. 

\begin{lemma} \label{mainit} Suppose $A\subseteq [1,N]$ with $|A|=\delta N$. If $$(A-A)\cap \left(I(g_{d_1})+ I(h_{d_2})\right)\subseteq \{0\}$$ and $d_1,d_{2},\delta^{-1}\leq \cQ$, then there exist $q\ll_{g,h} \delta^{-2}$ and $A'\subseteq [1,N']$ with 
$N'\gg_{g,h}  \delta^{2D(k\ell+1)}N$, 
\begin{equation*}\frac{|A'|}{N'} \geq (1+c\log^{-C}(\delta^{-1}))\delta   , \end{equation*} and $$(A'-A')\cap \left(I(g_{\lambda_2(q)d_1})+I(h_{\lambda_1(q)d_2})\right)\subseteq \{0\}, $$  for some $c=c(g,h)>0$ and $C=C(k,\ell)$. 
\end{lemma}

\subsection*{Proof of Theorem \ref{main2}} Throughout this proof, we let $C$ and $c$ denote sufficiently large or small positive constants, respectively, which we allow to change from line to line, but can depend only on $g$ and $h$. We use $C'$ and $c'$ similarly, but these constants can depend only on $k$ and $\ell$. 

\noindent Suppose $A \subseteq [1,N]$ with $|A|=\delta N$ and $$(A-A)\cap \left(I(g)+I(h)\right)\subseteq \{0\}.$$ Setting $A_0=A$, $N_0=N$, $d^{(0)}_1,d^{(0)}_2=1$, and $\delta_0=\delta$, Lemma \ref{mainit} yields, for each $m$, a set $A_m \subseteq [1,N_m]$ with $|A_m|=\delta_mN_m$ and $$(A_m-A_m)\cap \left(I\left(g_{d^{(m)}_1}\right)+I\left(h_{d^{(m)}_2}\right)\right)\subseteq \{0\}.$$ 

\noindent Further, we have that
\begin{equation} \label{NmI} N_m \geq c\delta^{2D(k\ell+1)}N_{m-1} \geq (c\delta)^{2D(k\ell+1)m} N,
\end{equation}

\begin{equation} \label{incsizeI} \delta_m \geq (1+c\log^{-C'}(\delta_{m-1}^{-1}))\delta_{m-1}, 
\end{equation}
and
\begin{equation}\label{dmI} d^{(m)}_i \leq (c\delta)^{-2k\ell} d^{(m-1)}_i \leq (c\delta)^{-2k\ell m},
\end{equation}
as long as 
\begin{equation} \label{delmI} d^{(m)}_i, \delta_m^{-1} \leq e^{\rho(\log N_m)^{1/3}} .
\end{equation}
However, we see that the density $\delta_m$ will exceed $1$, and hence (\ref{delmI}) must fail for $m=C\log^{C'}(\delta^{-1})$, which by (\ref{NmI}) and (\ref{dmI}) yields $(c\delta)^{-C\log^{C'}(\delta^{-1})}\geq e^{(\log N)^{1/3}},$ and hence $$\delta \ll_{g,h} e^{-(\log N)^{c'}}.$$This establishes Theorem \ref{main2} outside of the claim that we can take $c'=1/2$ if $\deg(g)=\deg(h)=2$, which we discuss in Section \ref{finalsec}. \qed

\subsection{Deducing Lemma \ref{mainit} from $L^2$ Fourier concentration} The philosophy behind the proof of Lemma \ref{mainit} is that a deficiency in polynomial differences in a set $A$ represents  nonrandom behavior, which should be detected in the Fourier analytic behavior of $A$. Specifically,  we locate one small denominator $q$ such that $\widehat{f_A}$ has $L^2$ concentration around rationals with denominator $q$, then use that information to find a long arithmetic progression on which $A$ has increased density. 

\begin{lemma}  \label{L2I} Suppose $A\subseteq [1,N]$ with $|A|=\delta N$, $\eta=c_0\delta$ for a sufficiently small constant $c_0=c_0(g,h)>0$, and  $\gamma=\eta^{-2D}/N$. If $(A-A)\cap \left(I(g_{d_1})+ I(h_{d_2})\right)\subseteq \{0\}$, $d_1,d_{2},\delta^{-1}\leq \cQ$, and  $|A\cap(N/9,8N/9)|\geq 3\delta N/4$, then there exists $q\leq \eta^{-2}$ such that 
\begin{equation*} \int_{\mathbf{M}'_q(\gamma)} |\widehat{f_A}(\alpha)|^2d\alpha \gg_{g,h} \delta^{2}\log^{-C}(\delta^{-1}) N  
\end{equation*}
for some $C=C(k,\ell)$.
\end{lemma}

\noindent Lemma \ref{mainit} follows from Lemma \ref{L2I} and the following standard $L^2$ density increment lemma, the continuous analog of Lemma \ref{dinczn}.

\begin{lemma}[Lemma 2.3 in \cite{thesis}, see also \cite{Lucier}, \cite{Ruz}] \label{dinc} Suppose $A \subseteq [1,N]$ with $|A|=\delta N$. If  $0< \theta \leq 1$, $q \in \N$, $\gamma>0$, and
\begin{equation*} \int_{\mathbf{M}'_q(\gamma)}|\widehat{f_A}(\alpha)|^2d\alpha \geq \theta\delta^2 N,
\end{equation*} 
then there exists an arithmetic progression 
\begin{equation*}P=\{x+\ell q : 1\leq \ell \leq L\}
\end{equation*}
with $qL \gg \min\{\theta N, \gamma^{-1}\}  $ and $|A\cap P| \geq \delta(1+\theta/32)L$.
\end{lemma}

\subsection*{Proof of Lemma \ref{mainit}}
 Suppose $A\subseteq [1,N]$, $|A|=\delta N$, $(A-A)\cap \left(I(g_{d_1})+ I(h_{d_2})\right)\subseteq \{0\}$, and $d_1,d_{2},\delta^{-1}\leq \cQ$. If $|A\cap (N/9,8N/9)| < 3\delta N/4$, then $\max \{ |A\cap[1,N/9]|, |A\cap [8N/9,N]| \} > \delta N/8$. In other words, $A$ has density at least $9\delta/8$ on one of these intervals. 

\noindent Otherwise, Lemmas \ref{L2I} and \ref{dinc} apply, so in either case, letting $\eta=c_0\delta$, there exists $q\leq \eta^{-2}$ and an arithmetic progression 
\begin{equation*}P=\{x+\ell q : 1\leq \ell \leq L\}
\end{equation*}
with $qL\gg_{g,h} \delta^{2D} N$ and $$|A\cap P|/L \geq (1+c\log^{-C}(\delta^{-1}))\delta  .$$ Partitioning $P$ into subprogressions of step size $\Lambda(q)=\lambda_1(\lambda_2(q))$, the pigeonhole principle yields a progression 
\begin{equation*} P'=\{y+a \Lambda(q) : 1\leq a \leq N'\} \subseteq P
\end{equation*}
with $N'\geq qL/2\Lambda(q)$ and $|A\cap P'|/N' \geq |A\cap P|/L$. This allows us to define a set $A' \subseteq [1,N']$ by \begin{equation*} A' = \{a \in [1,N'] : y+a \Lambda(q) \in A \},
\end{equation*} which satisfies $|A'|=|A\cap P'|$ and $N'\gg_{g,h} \delta^{2D}N/\Lambda(q) \gg_{g,h} \delta^{2D(k\ell+1)}N$. Moreover, by Proposition \ref{inh2}, $(A-A)\cap \left(I(g_{d_1})+ I(h_{d_2})\right)\subseteq \{0\}$ implies $(A'-A')\cap \left(I(g_{\lambda_2(q)d_1})+I(h_{\lambda_1(q)d_2})\right)\subseteq \{0\}$. \qed \\

Our task for this section is now completely reduced to a proof of Lemma \ref{L2I}.

\subsection{Preliminary notation for proof of Lemma \ref{L2I}} Before delving into the proof of Lemma \ref{L2I}, we take the opportunity to define some relevant sets and quantities, depending on our intersective polynomials $g,h\in \Z[x]$, scaling parameters $d_1,d_2$, a parameter $Y>0$, and the size of the ambient interval $N$. In all the notation defined below, we suppress all of the aforementioned dependence, as the relevant objects will be fixed in context.  

We define $W^{(1)}_{d_1}$, $\gamma_{d_1}^{(1)}$, and $j^{(1)}_{d_1}$ in terms of $g$ as in Section \ref{sievesec}. We then define $H_1$ to be the collection of natural number inputs $m\in W^{(1)}_{d_1}(Y)$ such that $g_{d_1}(m)$ is strictly between $0$ and $\pm N/18$, where the sign is the sign of the leading coefficient of $g$. We let $M_1=\left(\frac{N}{18|b|}\right)^{1/k}$,  where $b$ is the leading coefficient of $g_{d_1}$, and we let $$w_1=\prod_{p\leq Y}\left(1-\frac{j^{(1)}_{d_1}(p)}{p^{\gamma^{(1)}_{d_1}(p)}}\right). $$ 

We then analogously define $W^{(2)}_{d_2}$, $\gamma_{d_2}^{(2)}$, $j^{(2)}_{d_2}$, $H_2$, $M_2$, and $w_2$ in terms of $h$ and $d_2$, we let $$Z=\left\{(m,n)\in H_1\times H_2 : g_{d_1}(m)+h_{d_2}(n)=0\right\},$$ and we let $H=(H_1\times H_2)\setminus Z$. Letting $M=w_1w_2M_1M_2$, it follows from (\ref{symBIG}), (\ref{sieve}), and (\ref{logk}) that \begin{equation}\label{symdif5} \left|H_1 \triangle \left([1,M_1]\cap W^{(1)}_{d_1}(Y)\right)  \right| \ll_g 1, \end{equation} with the analogous statement for $H_2$, and \begin{equation}\label{Hsize} |H|\geq M/2, \end{equation} provided, for example, that $Y<e^{\sqrt{\log N}}$.

\subsection{Proof of Lemma \ref{L2I}} Suppose $A\subseteq [1,N]$ with $|A|=\delta N$, $(A-A)\cap \left(I(g_{d_1})+ I(h_{d_2})\right) \subseteq \{0\}$, and $d_1,d_2,\delta^{-1}\leq \cQ$.  Further, let $\eta=c_0\delta$ for an appropriately small $c_0=c_0(g,h)>0$, let $Q=\eta^{-2}$, and let $Y=\eta^{-2D}$.  Since $g_{d_1}(H_1)+h_{d_2}(H_2)\subseteq [-N/9,N/9]$, we have
\begin{align*} \sum_{\substack{x \in \Z \\ (m,n)\in H}} f_A(x)f_A(x+g_{d_1}(m)+h_{d_2}(n))&=\sum_{\substack{x \in \Z \\ (m,n)\in H}} 1_A(x)1_A(x+g_{d_1}(m)+h_{d_2}(n)) \\&-\delta\sum_{\substack{x \in \Z \\ (m,n)\in H}} 1_A(x)1_{[1,N]}(x+g_{d_1}(m)+h_{d_2}(n)) \\ &-\delta \sum_{\substack{x \in \Z \\ (m,n)\in H}} 1_{A}(x+g_{d_1}(m)+h_{d_2}(n))1_{[1,N]}(x) \\&+\delta^2\sum_{\substack{x \in \Z \\ (m,n)\in H}} 1_{[1,N]}(x)1_{[1,N]}(x+g_{d_1}(m)+h_{d_2}(n))  \\&\leq \Big(\delta^2N -2\delta|A\cap (N/9,8N/9)|\Big)|H|. 
\end{align*}
Therefore, if $|A \cap (N/9, 8N/9)| \geq 3\delta N/4$, then by (\ref{Hsize}) we have
\begin{equation}\label{neg} \sum_{\substack{x \in \Z \\ (m,n)\in H}} f_A(x)f_A(x+g_{d_1}(m)+h_{d_2}(n)) \leq -\delta^2NM/4.
\end{equation} 
We see from (\ref{symdif5}) and orthogonality of characters that 
\begin{equation}\label{orth}
\sum_{\substack{x \in \Z \\ (m,n)\in H}} f_A(x)f_A(x+g_{d_1}(m)+h_{d_2}(n))=\int_0^1 |\widehat{f_A}(\alpha)|^2 S(\alpha)d\alpha +O_{g,h}(N(w_1M_1+w_2M_2)),
\end{equation} 
where 
\begin{equation*}S_1(\alpha)= \sum_{\substack{m=1 \\ W^{(1)}_{d_1}(Y)}}^{M_1}e^{2\pi \i g_{d_1}(m)\alpha}, \quad S_2(\alpha)=\sum_{\substack{n=1 \\ W^{(2)}_{d_2}(Y)}}^{M_2}e^{2\pi \i h_{d_2}(n)\alpha}, \quad \text{and} \quad S(\alpha)=S_1(\alpha)S_2(\alpha).
\end{equation*} 
Combining (\ref{neg}) and (\ref{orth}), we have  
\begin{equation} \label{mass}
\int_0^1 |\widehat{f_A}(\alpha)|^2|S(\alpha)|d\alpha \geq \delta^2NM/8.
\end{equation} Letting $\gamma=\eta^{-2D}/N$, the estimates in Section \ref{expest} yield that if $d_1,d_2, \delta^{-1} \leq \cQ$, then for $\alpha \in \mathbf{M}_q(\gamma), \ q\leq Q $, we have \begin{equation} \label{SmajII} |S(\alpha)| \ll_{g,h} C^{\omega(q)}M/q,  
\end{equation} where $C=C(k,\ell)$. Further, for $\alpha \in \mathfrak{m}(\gamma,Q)$ we have 
\begin{equation} \label{SminII} |S(\alpha)| \leq  \delta M/16 ,
\end{equation} provided $c_0$ is chosen sufficiently small. The verification of (\ref{SmajII}) and (\ref{SminII}) is completely analogous to, though strictly easier than, the establishment of (\ref{Wmin}) and (\ref{Wmaj}) in Section \ref{it}.  

\noindent From (\ref{SminII}) and Plancherel's Identity, we have \begin{equation*}  \int_{\mathfrak{m}(\gamma,Q)} |\widehat{f_A}(\alpha)|^2|S(\alpha)|d\alpha \leq \delta^2NM/16, \end{equation*} which together with (\ref{mass}) yields \begin{equation}\label{majmass}  \int_{\mathfrak{M}(\gamma,Q)}|\widehat{f_A}(\alpha)|^2|S(\alpha)|d\alpha \geq \delta^2 NM/16. \end{equation} 
From (\ref{SmajII}) and (\ref{majmass}) , we have 
\begin{equation} \label{majmassII} \sum_{q=1}^Q  \frac{C^{\omega(q)}}{q}\int_{\mathbf{M}_q(\gamma)}|\widehat{f_A}(\alpha)|^2 {d}\alpha \gg_{g,h} \delta^2N.
\end{equation}
The function $b(q)=C^{\omega(q)}$ satisfies $b(qr)\geq b(r)$, and we make use of the following proposition.
\begin{proposition}\label{rstrick} For any $\gamma,Q>0$ satisfying $2\gamma Q^2<1$ and any function $b: \N \to [0,\infty)$ satisfying $b(qr)\geq b(r)$ for all $q,r\in \N$, we have $$\max_{q\leq Q} \int_{\mathbf{M}'_q(\gamma)}|\widehat{f_A}(\alpha)|^2 {d}\alpha \geq Q \Big(2\sum_{q=1}^Q b(q)\Big)^{-1} \sum_{r=1}^Q \frac{b(r)}{r}\int_{\mathbf{M}_r(\gamma)}|\widehat{f_A}(\alpha)|^2 {d}\alpha. $$
\end{proposition}

\begin{proof} By (\ref{majdisj}) we have \begin{align*} \Big(\sum_{q=1}^Q b(q)\Big) \max_{q\leq Q} \int_{\mathbf{M}'_q(\gamma)}|\widehat{f_A}(\alpha)|^2 {d}\alpha &\geq \sum_{q=1}^Q b(q) \int_{\mathbf{M}'_q(\gamma)}|\widehat{f_A}(\alpha)|^2 {d}\alpha \\ &= \sum_{q=1}^Q b(q) \sum_{r|q} \int_{\mathbf{M}_r(\gamma)}|\widehat{f_A}(\alpha)|^2 {d}\alpha \\ &= \sum_{r=1}^Q \int_{\mathbf{M}_r(\gamma)}|\widehat{f_A}(\alpha)|^2 {d}\alpha \sum_{q=1}^{Q/r} b(qr) \\ &\geq \frac{Q}{2}\sum_{r=1}^Q \frac{b(r)}{r}\int_{\mathbf{M}_r(\gamma)}|\widehat{f_A}(\alpha)|^2 {d}\alpha, 
\end{align*} 

\noindent where the last inequality comes from replacing $b(qr)$ with $b(r)$, and the proposition follows.
\end{proof}
\noindent Invoking the known estimate
\begin{equation*} \sum_{q=1}^Q C^{\omega(q)} \ll_C Q\log^C Q 
\end{equation*}
for any $C>0$ (see \cite{Selberg}), the lemma follows from (\ref{majmassII}) and Proposition \ref{rstrick}. \qed

\subsection{Discussion of the case $\deg(g)=\deg(h)=2$ and sums of $\ell \geq 3$ polynomials} \label{finalsec} In the case that $\deg(g)=\deg(h)=2$, the desired square root cancellation is already present in the complete exponential sums, and hence no sieving of inputs is required. This allows us to take $\cQ=\cQ(N)=N^{\rho}$ instead of $\cQ=e^{\rho(\log N)^{1/3}}$. Further, the $C^{\omega(q)}$ term is absent from the estimate (\ref{SmajII}), which allows us to replace $\log^{-C}(\delta^{-1})$ with a small constant $c$ in the conclusions of Lemmas \ref{mainit} and \ref{L2I}. Appropriately adjusting the proof in Section \ref{mainproof} yields $$\delta \ll_{g,h} e^{-c\sqrt{\log N}},$$ as claimed. 

\noindent If we consider sums of $\ell \geq 3$ intersective polynomials, then the square root cancellation in each variable allows us to replace (\ref{SmajII}) with $$|S(\alpha)| \ll_{h_1,\dots,h_{\ell}} C^{\omega(q)}M/q^{3/2}\ll M/q,$$ so we can again replace $\log^{-C}(\delta^{-1})$ with a small constant $c$ in the conclusions of Lemmas \ref{mainit} and \ref{L2I}. Further, if the reciprocals of the degrees of the polynomials add to at least 1, then, from the $q^{-1/\deg(h_i)}$ cancellation in the complete exponential sums in each variable, we get $$ |S(\alpha)| \ll_{h_1,\dots,h_{\ell}} M/q$$ without any sieving of inputs, so we can again take $\cQ=\cQ(N)=N^{\rho}$ instead of $\cQ=e^{\rho(\log N)^{1/3}}$. Appropriately adjusting the proof in Section \ref{mainproof} yields the following result, with which we conclude our discussion.

\begin{theorem} \label{main3} Suppose $\ell\geq 3$, $h_1,\dots, h_{\ell} \in \Z[x]$ are intersective polynomials, and $A\subseteq [1,N]$. If $$a-a'\neq \sum_{i=1}^{\ell}h_i(n_i)$$ for all distinct pairs $a,a'\in A$ and all $n_1,\dots, n_{\ell}\in \N$, then $$\frac{|A|}{N} \ll_{h_1,\dots,h_{\ell}} e^{-c(\log N)^{\mu}},$$ where $$\mu = \begin{cases}1/2 &\text{if } \sum_{i=1}^{\ell} \deg(h_i)^{-1} \geq 1 \\ 1/6 &\text{else} \end{cases}. $$
\end{theorem}

\noindent \textbf{Acknowledgements:} The author would like to thank the referee for their detailed and important recommendations,  Steve Gonek and Paul Pollack for their helpful comments and references, and Neil Lyall for his perpetual support.

\end{document}